\documentstyle{amltd}
\begin{document}
\annalsline{157}{2003}
 \received{May 29, 2001}
\startingpage{621}
\def\bye{\end{document}}
 \font\tenrm=cmr10
\input amssym.def
\input amssym.tex
\def\ritem#1{\item[{\rm #1}]}
\catcode`\@=11
\font\twelvemsb=msbm10 scaled 1100
\font\tenmsb=msbm10
\font\ninemsb=msbm10 scaled 800
\newfam\msbfam
\textfont\msbfam=\twelvemsb  \scriptfont\msbfam=\ninemsb
  \scriptscriptfont\msbfam=\ninemsb
\def\msb@{\hexnumber@\msbfam}
\def\Bbb{\relax\ifmmode\let\next\Bbb@\else
 \def\next{\errmessage{Use \string\Bbb\space only in math
mode}}\fi\next}
\def\Bbb@#1{{\Bbb@@{#1}}}
\def\Bbb@@#1{\fam\msbfam#1}
\catcode`\@=12

 \catcode`\@=11
\font\twelveeuf=eufm10 scaled 1100
\font\teneuf=eufm10
\font\nineeuf=eufm7 scaled 1100
\newfam\euffam
\textfont\euffam=\twelveeuf  \scriptfont\euffam=\teneuf
  \scriptscriptfont\euffam=\nineeuf
\def\euf@{\hexnumber@\euffam}
\def\frak{\relax\ifmmode\let\next\frak@\else
 \def\next{\errmessage{Use \string\frak\space only in math
mode}}\fi\next}
\def\frak@#1{{\frak@@{#1}}}
\def\frak@@#1{\fam\euffam#1}
\catcode`\@=12

\newcommand{\CC}{{\Bbb C}}
\newcommand{\RR}{{\Bbb R}}
\renewcommand{\SS}{{\Bbb S}}
\newcommand{\QQ}{{\Bbb Q}}
\newcommand{\FF}{{\Bbb F}}
\newcommand{\ZZ}{{\Bbb Z}}
\newcommand{\PP}{{\Bbb P}}
\renewcommand{\AA}{{\Bbb A}}
\newcommand{\GG}{{\Bbb G}}
\newcommand{\NN}{{\Bbb N}}
\newcommand{\HH}{{\Bbb H}}

\newcommand{\Zhat}{\hat{\ZZ}}
\newcommand{\AAf}{\AA_{\rm f}}
\newcommand{\Qbar}{{\overline{\QQ}}}
\newcommand{\Fbar}{{\overline{\FF}}}
\newcommand{\Spec}{{\rm Spec}}

\newcommand{\into}{\hookrightarrow}
\newcommand{\rH}{{\rm H}}
\newcommand{\calL}{{\cal L}}
\newcommand{\calV}{{\cal V}}
\newcommand{\calO}{{\cal O}}
\newcommand{\Supp}{{\rm Supp}}

\newcommand{\bs}{\backslash}
\newcommand{\ti}{\times}
\newcommand{\ol}{\overline}
\newcommand{\wh}{\widehat}
\newcommand{\wt}{\widetilde}
\newcommand{\lto}{\longrightarrow}
\newcommand{\lfrom}{\longleftarrow}
\newcommand{\sm}{{\rm sm}}

\newcommand{\End}{{\rm End}}
\newcommand{\Hom}{{\rm Hom}}
\newcommand{\Res}{{\rm Res}}
\newcommand{\Sh}{{\rm Sh}}
\newcommand{\Gal}{{\rm Gal}}
\newcommand{\ad}{{\rm ad}}
\newcommand{\ab}{{\rm ab}}
\newcommand{\der}{{\rm der}}
\newcommand{\nor}{{\rm nor}}
\newcommand{\inn}{{\rm inn}}
\newcommand{\GL}{{\rm GL}}
\newcommand{\bGL}{{\bf GL}}
\newcommand{\SL}{{\rm SL}}
\newcommand{\MT}{{\rm MT}}
\newcommand{\AM}{{\rm AM}}
\newcommand{\Gm}[2]{{\GG_{{\rm m}{{#1}}}^{{#2}}}}
\newcommand{\ld}{{\langle}}
\newcommand{\rd}{{\rangle}}
 
\relpenalty=10000
\binoppenalty=10000
 \title{Subvarieties of Shimura varieties}

 \acknowledgements{Both authors were
partially supported by the European Research Training Network Contract
HPRN-CT-2000-00120 ``arithmetic algebraic geometry''.}
 \twoauthors{Bas Edixhoven}{Andrei Yafaev}
\institutions{Mathematisch Instituut, Universiteit Leiden, The Netherlands\\
{\eightpoint {\it E-mail address\/}: edix@math.leidenuniv.nl}\\ \vglue2pt
Imperial College, London, England, UK\\
{\eightpoint {\it E-mail address\/}:  andrei.yafaev@ic.ac.uk}}
 
 \centerline{\it Dedicated to Laurent Moret-Bailly on the occasion of his $50^{\rm th}$ birthday}
 
\def\sni#1{\smallbreak\noindent{#1}. }
\def\ssni#1{\vglue-1pt\noindent\hskip18pt {#1}.}
\vfill
\centerline{\bf Contents}

\sni{1} Introduction
\sni{2} The strategy

\sni{3} Some preliminaries
\ssni{3.1} Mumford-Tate groups
\ssni{3.2} Variations of $\Bbb Z$-Hodge structure on Shimura varieties
\ssni{3.3} Representations of tori
\sni{4} Lower bounds for Galois orbits
\ssni{4.2} Galois orbits and Mumford-Tate groups
\ssni{4.3} Getting rid of $G$
\ssni{4.4} Proof of Proposition 4.3.9
\sni{5} Images under Hecke correspondences
\sni{6} Density of Hecke orbits
\sni{7} Proof of the main result
\ssni{7.3} The case where $i$ is bounded
\ssni{7.4} The case where $i$ is not bounded

\section{Introduction}\label{sec1}

The aim of this article is to prove a special case of the following
conjecture of Andr\'e and Oort on subvarieties of Shimura
varieties. For the terminology, notation, history and results obtained
so far we refer to the introduction of~\cite{Edixhoven1}, and the
references therein.
\eject

\proclaimtitle{Andr\'e-Oort}
\specialnumber{1.1} \proclaim{{C}onjecture} \label{AndreOortconjecture}
Let $(G,X)$ be a Shimura datum. Let $K$ be a compact open subgroup
of~$G(\AAf)$ and let $S$ be a set of special points in
$\Sh_K(G,X)(\CC)$. Then every irreducible component of the Zariski
closure of $S$ in $\Sh_K(G,X)_\CC$ is a subvariety of Hodge type.
\endproclaim

The choice of the special case that we will prove is motivated by work
of Wolfart~\cite{Wolfart1}  (see also Cohen and
W\"ustholz~\cite{CohenWustholz1}) on algebraicity of values of
hypergeometric functions at algebraic numbers. The hypergeometric
functions considered in \cite{Wolfart1} are the multi-valued
holomorphic $F(a,b,c)$ on $\PP^1(\CC)-\{0,1,\infty\}$ defined by:
$$
F(a,b,c)(z) = 1+\frac{ab}{c}z +
\frac{a(a+1)b(b+1)}{c(c+1)}\frac{z^2}{2!}+\cdots, \quad z\in\CC,\quad
|z|<1, 
$$
with $a$, $b$ and $c$ rational numbers, $-c$~not in~$\NN$. For the
following properties of the $F(a,b,c)$ the reader is referred
to~\cite{Wolfart1}. The functions $F(a,b,c)$ satisfy the differential
equations:
$$
z(z-1)F(a,b,c)''+((a+b+1)z-c)F(a,b,c)'+abF(a,b,c)=0. 
$$
Suppose from now on that~$a$, $b$, $c$, $a-c$ and $b-c$ are all not
integers.  Then, up to a factor in~$\Qbar^*$, $F(a,b,c)(z)$ is a quotient of a
   certain period $P(a,b,c)(z)$ of a certain abelian
   variety~$A(a,b,c,z)$, by a period $P_0(a,b,c)$ depending only
   on~$(a,b,c)$. This abelian variety $A(a,b,c,z)$ is a subvariety of
   the jacobian of the smooth projective model of the curve given by
   the equation $y^N=x^A(1-x)^B(1-zx)^C$, for suitable~$N$, $A$, $B$,
   and~$C$ depending on~$(a,b,c)$. More precisely, $P(a,b,c)(z)$ is the
integral of the differential form $y^{-1}dx$ (all of whose residues
are zero) over a suitable cycle. The fixed period $P_0(a,b,c)$ is a
period of an abelian variety of CM type. By \cite[Thm.~5]{Wustholz1} it follows that if $z$ and one value of
   $F(a,b,c)$ at $z$ are algebraic, then $A(a,b,c,z)$ is also of CM
   type, with the same type as~$P_0(a,b,c)$, and hence all values of
   $F(a,b,c)$ at $z$ are algebraic. (This theorem of W\"ustholz is
   about Grothendieck's conjecture on period relations and
   correspondences (see \cite{Andre1}) in the case of abelian
   varieties, and generalizes earlier work by Baker
   \cite[Cor.~2]{Wustholz1}.)
Because of this, it makes sense to ask the following question:
\begin{quote}
Under what conditions on $(a,b,c)$ is the set $E(a,b,c)$ of $z$ in\hfill\break
$\PP^1(\Qbar)-\{0,1,\infty\}$ such that $F(a,b,c)(z) \subset\Qbar$   finite?
\end{quote}
Wolfart proposes in \cite[Theorem]{Wolfart1} that the answer should
depend only on the monodromy group $\Delta(a,b,c)$ (with its
two-dimensional representation) of the differential equation. If this
monodromy group is finite, then $F(a,b,c)$ is algebraic over
$\Qbar(z)$, and hence $E(a,b,c)$ equals $\PP^1(\Qbar)-\{0,1,\infty\}$.
Suppose now that the monodromy group is infinite. Then one
distinguishes two cases: $\Delta(a,b,c)$ is arithmetic or not. In
terms of Shimura varieties, these two cases correspond to the image of
$\PP^1(\CC)-\{0,1,\infty\}$ in a suitable moduli space of polarized
abelian varieties under the map that sends $z$ to the isomorphism
class of $A(a,b,c,z)$ being of Hodge type or not
(see~\cite{CohenWustholz1}). If this image is of Hodge type, then the
set of $z$ in $\Qbar$ such that $A(a,b,c,z)$ is of CM type, with type  that of $P_0(a,b,c)$, is dense in $\PP^1(\CC)$ (even
for the Archimedean topology). Wolfart shows that under certain additional
conditions on $(a,b,c)$ the set $E(a,b,c)$ is infinite. Suppose now
that the image of $\PP^1(\CC)-\{0,1,\infty\}$ is not of Hodge
type. Then Wolfart's theorem claims that $E(a,b,c)$ is finite. But
Walter Gubler has pointed out an error in Wolfart's proof:
in~\cite[\S9]{Wolfart1}, there is no reason that the group
$\ol{\Delta}$ fixes the image of $\HH$ in the product of copies of the
unit disk under the product of the maps~$D_\omega$, and therefore the
identity $\delta^{\sigma}(f(w)) = f(\delta(w))$ for $\delta$ in
$\ol{\Delta}$ is not proved.

Since by W\"ustholz's theorem all $A(a,b,c,z)$ with $z$ in $E(a,b,c)$
are of a fixed CM type, hence isogeneous and hence contained in one
Hecke orbit, the following theorem completes Wolfart's program.

\vglue6pt {\elevensc Theorem 1.2.}
{\it Let $(G,X)$ be a Shimura datum and let $K$ be a compact open subgroup
of~$G(\AAf)$. Let $V$ be a finite\/{\rm -}\/dimensional faithful representation
of~$G${\rm ,} and for $h$ in $X$ let $V_h$ be the corresponding $\QQ$\/{\rm -}\/Hodge
structure. For $x=\ol{(h,g)}$ in $\Sh_K(G,X)(\CC)${\rm ,} let $[V_x]$ denote
the isomorphism class of~$V_h$. Let $Z$ be an irreducible closed
algebraic curve contained in $\Sh_K(G,X)_\CC$ such that $Z(\CC)$
contains an infinite set of special points $x$ such that all $[V_x]$
are equal. Then $Z$ is of Hodge type. In particular{\rm ,} if $Z$ is an
irreducible closed algebraic curve contained in $\Sh_K(G,X)_\CC$ such
that $Z(\CC)$ contains an infinite set of special points that lie in
one Hecke orbit{\rm ,} then $Z$ is of Hodge type.}
\vglue6pt
\advance\theoremcount by 2

This theorem, in the case where the special points in question lie in
one Hecke orbit, was first proved in the second author's
thesis~\cite{Yafaev1}, in which one chapter (providing a lower bound
for Galois orbits) was written by the first author. The main
difference between this article and the thesis is that now we consider
isomorphism classes of $\QQ$-Hodge structures instead of Hecke
orbits. This makes it possible to reduce the proof of the theorem to
the case where $Z$ is Hodge generic and $G$ of adjoint type (the proof
in the thesis could not achieve this and was therefore more difficult
to follow). 

The proof given in this article is nice because it is entirely in
``$(G,X)$-language''; the main tools are algebraic groups and their
groups of adelic points.  But it is not completely satisfactory in the
sense that it should be possible to proceed as in~\cite{Edixhoven1},
i.e., without distinguishing the two cases as we do in
Section~\ref{sec6}. On the other hand, the proof in the first of these
two cases can lead to a generalization to arbitrary Shimura varieties
of Moonen's result in~\cite[\S5]{Moonen3}
(Conjecture~\ref{AndreOortconjecture} for moduli spaces of abelian
varieties, and sets of special points for which there is a prime at
which they are all ``canonical''). Finally, it would be nice to
replace the condition that all $[V_x]$ are the same in the theorem
above by the condition that all associated Mumford-Tate groups are
isomorphic (this would give a statement that does not depend on the
choice of a representation).
 
 \vglue-12pt
\section{The strategy}\label{sec2}
 \vglue-6pt
The aim of this section is to explain the strategy of the proof of
Theorem~1.2. In Sections~\ref{sec3}--\ref{sec5} we will
prove the necessary ingredients to  be put together in
Section~\ref{sec6}. 

We observe that the compact open subgroup $K$ in
Conjecture~\ref{AndreOortconjecture} is irrelevant: for $(G,X)$ a
Shimura datum, $K$ and $K'$ open compact subgroups of $G(\AAf)$ with
$K\subset K'$, an irreducible subvariety $Z$ of $\Sh_K(G,X)_\CC$ is of
Hodge type if and only if its image in $\Sh_{K'}(G,X)_\CC$ is. A bit
more generally, for $(G,X)$ a Shimura datum, $K$ and $K'$ open compact
subgroups of~$G(\AAf)$ and $g$ in~$G(\AAf)$, an irreducible subvariety
$Z$ of $\Sh_K(G,X)_\CC$ is of Hodge type if and only if one (or
equivalently, all) of the irreducible components of $T_gZ$ is (are) of
Hodge type, where $T_g$ is the correspondence from $\Sh_K(G,X)_\CC$ to
$\Sh_{K'}(G,X)_\CC$ induced by~$g$. The
irreducible components of intersections of subvarieties of Hodge type
are again of Hodge type (this is clear from the interpretation of
subvarieties of Hodge type as loci where certain classes are Hodge
classes). Hence there does exist a smallest subvariety of Hodge type
of $\Sh_K(G,X)_\CC$ that contains $Z$; our first concern is now to
describe that subvariety.

\proclaim{Proposition}\label{prop2.1}
Let $(G,X)$ be a Shimura datum{\rm ,} $K$ a compact open subgroup of
$G(\AAf)$ and let $Z$ be a closed irreducible subvariety of\/
$\Sh_K(G,X)_\CC$. Let $s$ be a Hodge generic point of~$Z${\rm :} its
Mumford\/{\rm -}\/Tate group is the generic Mumford\/{\rm -}\/Tate group on~$Z$. Let
$(x,g)$ in $X\ti G(\AAf)$ lie over $s${\rm ,} and let $G'$ be the
Mumford\/{\rm -}\/Tate group of~$x$. Then we have a morphism of Shimura data
from $(G',X')$ to $(G,X)$ with $X'$ the $G'(\RR)$\/{\rm -}\/conjugacy class
of~$x$. Let $K'$ be the intersection of $G'(\AAf)$ and
$gKg^{-1}$. Then the inclusion of $G'$ in $G${\rm ,} followed by right
multiplication by $g$ induces a morphism
$f\colon\Sh_{K'}(G',X')_\CC\to\Sh_K(G,X)_\CC$. This morphism is finite
and its image contains~$Z$. Let $Z'$ be an irreducible component
of~$f^{-1}Z$. Then $Z$ is of Hodge type if and only if $Z'$ is.
\endproclaim

{\it Proof}.
This follows from Proposition~2.8 and Section~2.9 of \cite{Moonen2}.
\hfill\qed\vglue4pt

Proposition~\ref{prop2.1} shows that
Conjecture~\ref{AndreOortconjecture} is true if and only if it is true
for all sets of special points $S$ whose Zariski closure is
irreducible and Hodge generic. Similarly, Proposition~\ref{prop2.1}
reduces the proof of Theorem~1.2 to the case where $Z$ is
Hodge generic. We note that even if $Z(\CC)$ has an infinite
intersection with the Hecke orbit of a special point, this is not
necessarily so for~$Z'$, because the inverse image in
$\Sh_{K'}(G',X')_\CC$ of a Hecke orbit in $\Sh_K(G,X)_\CC$ is a
disjoint union of a possibly infinite number of Hecke orbits. This
explains why we work with equivalence classes of $\QQ$-Hodge
structures.

\proclaim{Proposition}\label{prop2.2}
Let $(G,X)$ be a Shimura datum{\rm ,} let $G^\ad$ be the quotient of $G$ by
its center{\rm ,} and let $X^\ad$ be the $G^\ad(\RR)$\/{\rm -}\/conjugacy class of
morphisms from $\SS$ to $G^\ad_\RR$ that contains the image
of~$X$. Let $K^\ad$ be a compact open subgroup of $G^\ad(\AAf^{\phantom{|}})${\rm ,} and
let $K$ be a compact open subgroup of $G(\AAf)$ whose image in
$G^\ad(\AAf)$ is contained in~$K^\ad$. Then the induced morphism from
$\Sh_K(G,X)_\CC$ to $\Sh_{K^\ad}(G^\ad,X^\ad)_\CC$ is finite. Let $Z$
be a closed irreducible subvariety of\/$\Sh_K(G,X)_\CC${\rm ,} and let
$Z^\ad$ be its image in\/$\Sh_{K^\ad}(G^\ad,X^\ad)_\CC$. Then $Z$ is of
Hodge type if and only if $Z^\ad$ is. 
\endproclaim

{\it Proof}.
By \cite[\S2.1]{Moonen2}, $X$ is just a union of connected components
of~$X^\ad$. Let $S$ and $S^\ad$ be the connected components of $\Sh_K(G,X)_\CC$ and of\break
$\Sh_{K^\ad}(G^\ad,X^\ad)_\CC$ that contain
$Z$ and $Z^\ad$, respectively. Let $X^+$ be a connected component of
$X$ and let $g$ in $G(\AAf)$ be such that $S$
is the image in $\Sh_K(G,X)_\CC$ of $X^+\times\{g\}$. Then the inverse
images of $Z$ and $Z^\ad$ in $X^+\times\{g\}$ are equal, hence the
property of being of Hodge type for them is equivalent. 
\hfill\qed
\vglue4pt

We want to use Proposition~\ref{prop2.2} to reduce the proof of
Theorem~1.2 to the case where $G$ is semi-simple of adjoint
type. In order to do that, all we need to do is to construct a
faithful representation $W$ of $G^\ad$ such that $Z^\ad(\CC)$ contains
a Zariski dense set of special points $x$ with all $[W_x]$ equal.

\numbereddemo{Construction}\label{cons2.3}
Let $G$ be a reductive algebraic group over $\QQ$ and let $V$ be a
faithful finite-dimensional representation. Let $E$ be a finite
extension of $\QQ$ such that the representation $V_E$ of $G_E$ is a
direct sum $V_1\oplus\cdots\oplus V_r$ with each $V_i$ absolutely
irreducible. Let $C$ be the center of~$G$; then $C_E$ acts via a
character $\chi_i$ on~$V_i$. For each~$i$, we let $d_i=\dim_E(V_i)$, and
we define $W_i:=V_i^{\otimes_E d_i}\otimes_E\det_E(V_i)^*$. Then
$W:=W_1\oplus\cdots\oplus W_r$ is a faithful representation
of~$G_E^\ad$. We get a faithful representation of $G^\ad$ on $W$ as
$\QQ$-vector space via the sequence of injective morphisms of
algebraic groups: 
\smallbreak
\hfil  ${\displaystyle
G^\ad\to\Res_{E/\QQ}G_E^\ad\to\Res_{E/\QQ}{\bf GL}_E(W)\to{\bf GL}_\QQ(W).
}$\hfill
\enddemo

Suppose now that $\Sigma$ is a Zariski dense subset of
$\Sh_K(G,X)_\CC$ such that the $\QQ$-Hodge structures $V_x$ with $x$
in $\Sigma$ are all isomorphic to a fixed $\QQ$-Hodge
structure~$H$. Then the $\QQ$-Hodge structures with $E$-coefficients
$E\otimes V_x$ are all isomorphic to~$E\otimes H$. Now $E\otimes H$ is
a direct sum of finitely many simple $\QQ$-Hodge structures with
$E$-coefficients. Hence there are, up to isomorphism, only finitely
many ways to decompose $E\otimes H$ into a direct sum of $r$
terms. Hence $\Sigma$ is a finite disjoint union of subsets $\Sigma_i$
such that for each $i$ the $E\otimes V_x=V_{1,x}\oplus\cdots\oplus
V_{r,x}$ with $x$ in $\Sigma_i$ are all isomorphic term by term. It
follows that the $W_x$ with $x$ in $\Sigma_i$ are all isomorphic. So
it remains to prove Theorem~1.2 in the case where $G$ is
semi-simple of adjoint type and $Z$ Hodge generic.

At this point we can describe the strategy of the proof of the main
result. So let the notation now be as in Theorem~1.2, and
suppose that $G$ is semi-simple of adjoint type, that $Z$ is Hodge
generic, and that $K$ is neat. Let $\Sigma$ be an infinite subset of
special points of $Z$ such that all $[V_x]$ with $x$ in $\Sigma$ are
equal. Let $S$ be the irreducible component of $\Sh_K(G,X)$ that
contains~$Z$. We will show that there exists $g$ in $G(\AAf)$ such that
an irreducible component $T_g^0$ of the Hecke correspondence on $S$
induced by $g$ has the property that $T_g^0Z=Z=T_{g^{-1}}^{0^{\phantom{1}}}Z$ (with
$T_{g^{-1}}^0$ the transpose of~$T_g^0$) and is such that all the
$T_g^0+T_{g^{-1}}^{0^{\phantom{1}}}$-orbits in $S$ are dense (for the Archimedean
topology). This clearly implies that $Z=S$, so that $Z$ is of Hodge
type.

To find such a~$g$, we proceed as follows. We will take $g$ in
$G(\QQ_p)$, for some prime number~$p$. For all but finitely many $p$,
the image of $Z$ under each irreducible component of any $T_g$ with
$g$ in $G(\QQ_p)$ is either empty or irreducible. The proof of this
will be given in Section~\ref{sec4}, whose main ingredient is
Theorem~\ref{thmNori} by Nori. The density of all
$T_g^0+T_{g^{-1}_{\phantom{1}}}^0$-orbits will be proved in Section~\ref{sec5},
under the assumption that no image of $g$ under projection to a simple
factor $H$ of $G$ is contained in a compact subgroup of~$H(\QQ_p)$. To
get the equalities $T_g^0 Z=Z$ and $T_{g^{-1}}^0 Z=Z$ we try to find
$g$ such that $T_gZ\cap Z$ contains a large number of the given
special points, compared to the degree of the correspondence~$T_g$. In
doing this, we distinguish two cases. In one case, the intersection
will contain at least one big Galois orbit. In the other case, it
contains infinitely many of the given special points. The main
ingredient here is the description of the Galois action on special
points, plus a lower bound on the number of points in the Galois
orbits of our given special points that will be established in
Section~\ref{sec3}.

\vglue-6pt
\section{Some preliminaries}
\vglue-6pt

 3.1. {\it Mumford\/{\rm -}\/Tate groups}.
For $V$ a free $\ZZ$-module of finite rank, we define $\bGL(V)$ to be
the group scheme given by $\bGL(V)(A)=\GL_A(V_A)$ for all
rings~$A$. For $h$ a $\ZZ$-Hodge structure, i.e., a free $\ZZ$-module
of finite rank, together with a morphism $h\colon\SS\to\bGL(V)_\RR$,
we let $\MT(V,h)$ be the Zariski closure in $\bGL(V)$ of the usual
Mumford-Tate group $\MT(V_\QQ,h)$ in $\bGL(V)_\QQ$.

\demo{{\rm 3.2.} Variations of $\ZZ$-Hodge structure on Shimura varieties}
Let $(G,X)$ be a Shimura datum, $K$ a neat compact open subgroup of
$G(\AAf)$ and $\rho\colon G\to\GL_n$ a representation that factors
through $G\to G^\ad$, such that $\rho(K)$ is in $\GL_n(\Zhat)$. Then
there is a variation of $\ZZ$-Hodge structures $V$ on $\Sh_K(G,X)$
constructed as follows. On $X\times G(\AAf)/K$, we consider the
variation of $\ZZ$-Hodge structure $V_1$ whose restriction to
$X\times\{\ol{g}\}$ is $\QQ^n\cap\rho(g)(\Zhat^n)\times X$ (with the
$\QQ$-Hodge structure on $\QQ^n\times\{x\}$   given by the morphism
$\rho_\RR\circ x$ from $\SS$ to~$\GL_{n,\RR}$). Then $G(\QQ)$ acts on
$V_1$, and the quotient is the $V$ that we want (for each $(x,\ol{g})$
in $X\times G(\AAf)/K$ and $q$ in $G(\QQ)$ that stabilizes
$(x,\ol{g})$, the image of $q$ in $G^\ad(\QQ)$ is trivial, hence
$\rho(q)$ is the identity).

A more conceptual way to describe $V$ is as follows. We consider two
actions of $G(\QQ)\times G(\AAf)$ on $\AAf^n\times X\times G(\AAf)$
given by:
\begin{eqnarray*}
(q,k)*_1(v,x,g) & = & (qv,qx,qgk), \\
(q,k)*_2(v,x,g) & = & (k^{-1}v,qx,qgk).
\end{eqnarray*}
The first action stabilizes $\QQ^n\times X\times G(\AAf)$, and the
second stabilizes $\Zhat^n\times X\break\times G(\AAf)$. The quotient by the
first action gives a locally constant sheaf $V_\QQ$ of $\QQ$-vector
spaces on $\Sh_K(G,X)(\CC)$, and the second one a locally constant
sheaf $V_{\Zhat}$ of $\Zhat$-modules.  The automorphism:
$$
\AAf^n\times X\times G(\AAf) \lto \AAf^n\times X\times G(\AAf), \quad
(v,x,g)\mapsto (g^{-1}v,x,g)
$$
transforms the first action into the second, hence gives an
isomorphism between the two locally constant sheaves of $\AAf$-modules
on $\Sh_K(G,X)(\CC)$. Then $V$ is the ``intersection'' of $V_\QQ$
and $V_{\Zhat}$ in~$V_{\AAf}$, i.e., the inverse image under this
isomorphism of $V_{\Zhat}$ in~$V_\QQ$.
\enddemo

3.3. {\it Representations of tori.} 
A torus over a scheme $S$ is an $S$-group scheme $T$ that is of the
form $\Gm{S}{r}$, locally for the fpqc topology on~$S$ (\cite[Exp.~IX,
D\'ef.~1.3]{Grothendieck1}). If $S$ is normal and noetherian, then a
torus $T/S$ is split over a suitable surjective finite \'etale cover of
$S'\to S$; i.e., $T_{S'}$ is isomorphic to some $\Gm{S'}{r}$
(\cite[Exp.~IX, Thm.~5.16]{Grothendieck1}); one may take $S'\to S$
Galois and $S'$ connected. If $S$ is integral normal and noetherian,
with generic point~$\eta$, then any isomorphism $f\colon T_{1,\eta}\to
T_{2,\eta}$ with $T_1$ and $T_2$ tori over $S$ extends uniquely to an
isomorphism over~$S$ (use~\cite[Exp.~X, Cor.~1.2]{Grothendieck1}).

For $S$ a connected scheme,
$T\cong\Gm{S}{r}$ a split torus and $V$ an $\calO_S$-module,\break it is
equivalent to give a $T$-action on $V$ or an $X$-grading on $V$, with
$X=\break\Hom(T,\Gm{S}{})$ the character group of~$T$ (\cite[Exp.~I,
\S4.7]{Grothendieck1}).

Suppose now that $S$ is an integral normal noetherian scheme, that $T$
is a torus over $S$ and that $\pi\colon S'\to S$ is a connected finite \'etale Galois cover with group~$\Gamma$ over which $T$ is split. Let
$X$ be the character group of~$T_{S'}$; then $X$ is a free
$\ZZ$-module of finite rank with a $\Gamma$-action. Then, to give an
action of $T$ on a quasi-coherent $\calO_S$-module $V$ is equivalent
to giving an $X$-grading $V_{S'}=\pi^* V=\oplus_{x} V_{S',x}$ such
that for all $\gamma$ in $\Gamma$ and all $x$ in $X$ one has $\gamma
V_{S',x}=V_{S',\gamma x}$ (to see this, use finite \'etale descent of
quasi-coherent modules as in \cite[\S6.2]{BLR1}). If $V$ and $W$ are
two representations of $T$ on locally free $\calO_S$-modules of finite
rank, then $V$ and $W$ are isomorphic, locally for the Zariski
topology on $S$, if and only if for all $x$ in $X$ the ranks of
$V_{S',x}$ and $W_{S',x}$ are equal (use that $\Hom_{\calO_S}(V,W)^T$
is a direct summand of $\Hom_{\calO_S}(V,W)$, whose formation commutes
with base change).

\specialnumber{3.3.1}\proclaim{Lemma}\label{lem:condatp}
Let $p$ be a prime number{\rm ,} let $T$ be a torus over~$\ZZ_p${\rm ,} let $V$ be
a free $\ZZ_p$\/{\rm -}\/module of finite rank equipped with a faithful action
of $T_{\QQ_p}$ on~$V_{\QQ_p}$. Let $T'$ be the scheme-theoretic
closure of $T_{\QQ_p}$ in~$\bGL(V)$. Then the following conditions are
equivalent\/{\rm :}\/
\begin{itemize}
\ritem{1.} $T'_{\FF_p}$ is a torus\/{\rm ;}\/
\ritem{2.} $T'$ is a torus\/{\rm ;}\/
\ritem{3.} The action of\/ $T_{\QQ_p}$ on $V_{\QQ_p}$ extends to an action
of\/ $T$ on~$V$\/{\rm ;}\/
\ritem{4.} $T$ stabilizes the lattice $V$ in the sense that for all finite
extensions $K$ of\/ $\QQ_p$ the lattice $V_{O_K}$ is stabilized by all
elements of\/~$T(O_K)$.
\end{itemize}
The set of $\ZZ_p$-lattices in $V_{\QQ_p}$ that are fixed by $T$ form
exactly one orbit under $C(\QQ_p)${\rm ,} where $C$ denotes the
centralizer of $T$ in~$\bGL(V)$.
\endproclaim
\demo{Proof}
Suppose first that $T'_{\FF_p}$ is a torus. Then $T'$, being a flat
group scheme affine and of finite type over $\ZZ_p$ whose fibers over
$\FF_p$ and $\QQ_p$ are tori, is a torus by \cite[Exp.~X,
Cor.~4.9]{Grothendieck1}.

Suppose that $T'$ is a torus. Then $T'=T$ by~\cite[Exp.~X,
Cor.~1.2]{Grothendieck1}). Hence the action of $T_{\QQ_p}$ on
$V_{\QQ_p}$ extends to an action of $T$ on~$V$.

Now suppose that the action of $T_{\QQ_p}$ on $V_{\QQ_p}$ extends to
an action of $T$ on~$V$. Then $T$ stabilizes~$V$, by definition. Also,
the description above of representations of tori shows that $T$ acts
faithfully on~$V$, so that $T$ is a closed subscheme of $\bGL(V)$,
flat over $\ZZ_p$, and hence equal to the scheme-theoretic closure of
its generic fiber. So $T'=T$ and $T'_{\FF_p}$ is a torus.

Suppose that $T$ stabilizes~$V$. Let $K$ be the splitting
field of~$T_{\QQ_p}$. Then $T_{O_K}$ is a split torus, and the action
of $T_{\QQ_p}$ on $V_{\QQ_p}$ is given by an $X$-grading of~$V_K$,
where $X$ is the character group of~$T_K$. Let $m$ be an integer that
is prime to $p$ such that the characters $x$ in $X$ with $V_{K,x}\neq
0$ have distinct images in~$X/mX$. Since $T$ stabilizes~$V$, the
$m$-torsion subgroup scheme $T[m]$ of $T$ acts on~$V$. This action
corresponds to an $X/mX$-grading on $V_{O_K}$ that is compatible with
the $X$-grading on~$V_K$. Hence the $X$-grading on $V_K$ extends to an
$X$-grading on $V_{O_K}$, which shows that the action of $T_{\QQ_p}$
on $V_{\QQ_p}$ extends to an action of $T$ on~$V$. 

Finally, let $S$ be the set of $\ZZ_p$-lattices in $V_{\QQ_p}$ that
are fixed by~$T$. Let $W$ be any $\ZZ_p$-lattice in~$V_{\QQ_p}$. The
$T_{\QQ_p}$-action on $W_{\QQ_p}=V_{\QQ_p}$ corresponds to the
$X$-grading on~$V_K$. By finite \'etale descent, the $O_K$-submodule
$\oplus_x (W_{O_K}\cap V_{K,x})$ of $W_{O_K}$ is of the form
$W'_{O_K}$ for a unique $\ZZ_p$-lattice $W'$ contained in~$W$. Then
$W'_{O_K}$ is the direct sum of the $W'_{O_K}\cap V_{K,x}$, hence is a
representation of~$T$. Hence $W'$ is in~$S$; in fact, $W'$ is the
largest sublattice of $W$ that is fixed by~$T$. In particular, $S$ is
not empty. Let now $V_1$ and $V_2$ be two \pagebreak elements of~$S$. Then both
are representations of~$T$. Since for each $x$ the $V_{i,O_K,x}$ are
of equal rank, $V_1$ and $V_2$ are isomorphic as representations
of~$T$. Let $g\colon V_1\to V_2$ be an isomorphism. Then $g$ is an
element of $C(\QQ_p)$ that sends $V_1$ to~$V_2$.
\enddemo

\vglue-12pt
\section{Lower bounds for Galois orbits}\label{sec3}

The aim of this section is to give certain lower bounds for the sizes
of Galois orbits of special points on a Shimura variety. To be
precise, we will prove the following theorem.

\proclaim{Theorem}\label{thm3.1}
Let $(G,X)$ be a Shimura datum{\rm ,} with $G$ semi\/{\rm -}\/simple of adjoint type{\rm ,}
and let $K$ be a neat compact open subgroup of\/~$G(\AAf)$. Via a
suitable faithful representation we view $G$ as a closed algebraic
subgroup of\/~$\GL_{n,\QQ}${\rm ,} such that $K$ is in~$\GL_n(\Zhat)$. Let
$V$ be the induced variation of $\ZZ$\/{\rm -}\/Hodge structures
on\/~$\Sh_K(G,X)_\CC$.  Let $s_0$ be a special point
of\/~$\Sh_K(G,X)_\CC$. Let $F\subset\CC$ be a number field over which
the Shimura variety\/ $\Sh_K(G,X)_\CC$ has a canonical model\/
$\Sh_K(G,X)_F${\rm ;} i.e.{\rm ,} a finite extension of the reflex field
associated to~$(G,X)$. Then there exist real numbers $c_1>0$ and
$c_2>0$ such that for all $s$ in $\Sh_K(G,X)_F(\Qbar)$ such that the
$\QQ$\/{\rm -}\/Hodge structure $V_{s,\QQ}$ is isomorphic to $V_{s_0,\QQ}$, we
have\/{\rm :}\/
$$
|\Gal(\Qbar/F){\cdot}s| > c_1\prod_{\{\hbox{$p$
prime}\;|\;\hbox{\it $\MT(V_s)_{\FF_p}$ is not a torus}\}} c_2p.
$$
\endproclaim 

Let us note that varying $F$ and $K$ does not affect the statement of
the theorem: if $F'$ and $K'$ satisfy the same hypotheses as $F$ and
$K$, then the sizes of the Galois orbits differ by a bounded factor,
and $K'\cap K$ has finite index in both $K$ and~$K'$. In the course of
the proof of Theorem~\ref{thm3.1} we will assume that $F$ is the
splitting field of $M_\QQ$, with $M$ the Mumford-Tate group of~$s_0$.

We note that $M(\RR)$ is compact: the kernel of the action of~$G(\RR)$
on $X$ consists precisely of the product of the compact factors, i.e.,
if $G_\RR$ is the product of simple $G_i$'s, then the kernel is the
product of the $G_i(\RR)$ that are compact, and $M(\RR)$ stabilizes a
point in the hermitian manifold~$X$. It follows that $M(\QQ)$ is
discrete in~$M(\AAf)$.

\demo{{\rm 4.2.} Galois orbits and Mumford-Tate groups} 
Let the notation be as in Theorem~\ref{thm3.1}. We choose a set of
representatives $R$ in $G(\AAf)$ for the quotient $G(\QQ)\backslash
G(\AAf)/K$; note that $R$ is finite. Then for $s$ in
$\Sh_K(G,X)_F(\CC)$ there exists a unique $g_s$ in $R$ and an element
$\tilde{s}$ in $X$ unique up to $\Gamma_s:=G(\QQ)\cap g_sKg_s^{-1}$,
such that $s=\ol{(\tilde{s},g_s)}$. We fix a choice for~$\tilde{s}_0$.
\vglue3pt

Let $s=\ol{(\tilde{s},g_s)}$ be in $\Sh_K(G,X)_F(\Qbar)$ such that the
$\QQ$-Hodge structure $V_{s,\QQ}$ is isomorphic to~$V_{s_0,\QQ}$. Then
$\tilde{s}$ gives an embedding of the Mumford-Tate group $\MT(s)_\QQ$
in $G$, and an inclusion of Shimura data from $(\MT(s)_\QQ,\{\tilde{s}\})$
in~$(G,X)$. Note that $\MT(s)_\QQ$ is isomorphic to~$M_\QQ$, hence has
splitting field~$F$. This gives morphisms of Shimura varieties
over~$F$:
\begin{equation}\label{eqn3.2.1}
\Sh_{K_s\cap\MT(s)(\AAf)}(\MT(s)_\QQ)_F \lto \Sh_{K_s}(G,X)_F
\stackrel{{\cdot}g_s}{\lto} \Sh_K(G,X)_F, \qquad \speqnu{4.2.1}
\end{equation} 
with the first one given by the inclusion, the second one by right
multiplication by~$g_s$, and with $K_s=K\cap g_sKg_s^{-1}$. By
construction, $V_{s,\QQ}$ and $V_{s_0,\QQ}$ are both $\QQ^n$ with
Hodge structures $\tilde{s}$ and $\tilde{s}_0\colon\SS\to\GL_{n,\RR}$,
respectively. Let $h$ in $\GL_n(\QQ)$ be an automorphism of $\QQ^n$
that is an isomorphism from $V_{s_0,\QQ}$ to~$V_{s,\QQ}$. Then we have
$\tilde{s}=\inn_h\circ\tilde{s}_0$. It follows that the reciprocity
morphisms $r_{\tilde{s}}$ and~$r_{\tilde{s}_0}$, viewed as morphisms
of tori over $\QQ$ from $\Res_{F/\QQ}\Gm{F}{}$ to $G$ with images
$\MT(s)_\QQ$ and $M_\QQ$ are related by $r_{\tilde{s}}=\inn_h\circ
r_{\tilde{s}_0}$. In particular, the isomorphism from $M_\QQ$ to
$\MT(s)_\QQ$ induced by $\inn_h$ gives an isomorphism of Shimura
varieties over~$F$:
\begin{equation}
\Sh(M_\QQ)_F \lto \Sh(\MT(s)_\QQ)_F. \speqnu{4.2.2} \label{eqn:isomSHMT}
\end{equation} 
The Galois group $\Gal(\Qbar/F)$ acts on $\Sh(M_\QQ)=M(\QQ)\backslash
M(\AAf)$ via its maximal abelian quotient, which we view, via class
field theory, as a quotient of $(\AA\otimes F)^*/(\RR\otimes
F)^{*,+}$. The action of $\Gal(\Qbar/F)$ is then given by
$r_{\tilde{s}_0}\colon \Res_{F/\QQ}\Gm{F}{}\break\to M_\QQ$.  We note that
$(\AA\otimes F)^*/(\RR\otimes F)^{*,+}$ is the product of
$(\AAf\otimes F)^*$ and the finite group of connected components of
$(\RR\otimes F)^*$. Hence the size of the $\Gal(\Qbar/F)$-orbit of $s$
in the $\Sh_K(G,X)_F(\Qbar)$ is, up to a bounded factor which is
independent of~$s$, the size of the $(\AAf\otimes F)^*$-orbit. Since
$F^{*,+}$ in $(\AAf\otimes F)^*$ acts trivially, and since the class
group $F^*\backslash(\AAf\otimes F)^*/(\Zhat\otimes O_F)^*$ is finite,
proving the lower bound we want for the $(\AAf\otimes F)^*$-orbits is
equivalent to proving it for the $(\Zhat\otimes
O_F)^*$-orbits. Moreover, since the set $R$ is finite, and since each
${{\cdot}g_s}\colon\Sh_{K_s}(G,X)_F\to\Sh_K(G,X)_F$ is finite, it is
enough to prove the lower bound for the $(\Zhat\otimes O_F)^*$-orbit
in $\Sh_{K_s}(G,X)_F(\Qbar)$.
\enddemo

4.3. {\it Getting rid of $G$.}
In order to simplify our task (i.e., to prove Theorem~\ref{thm3.1}) we
introduce the pair $(\GL_{n,\QQ},Y)$, with $Y$ the $\GL_n(\RR)$
conjugacy class in $\Hom_\RR(\SS,\GL_{n,\RR})$ that contains the image
of $X$ under~$\rho$. Of course, this pair is not a Shimura datum (if
$n>2$). 

For $s=\ol{(\tilde{s},g_s)}$ in $\Sh_K(G,X)_F(\Qbar)$ and $h$ in
$\GL_n(\QQ)$ as in the previous section, we consider the following
commutative diagram (of sets):
$$
\begin{array}{ccl}
\Sh_{K_s\cap\MT(s)(\AAf)}(\MT(s)_\QQ)(\CC) &\longrightarrow&
 \Sh_{K_s}(G,X)(\CC)
  \\[4pt]
{\scriptstyle i_s}\Big\uparrow&&\hskip.45in\Big\downarrow\\ [4pt]
\Sh_{K_s'}(M_\QQ)(\CC) &\stackrel{f_s}{\longrightarrow}&
\GL_n(\QQ)\backslash\left(Y\times\GL_n(\AAf)\right)/\GL_n(\Zhat),
\end{array}
$$  
where $K_s'$ is the subgroup of $M(\AAf)$ that corresponds to
$K_s\cap\MT(s)(\AAf)$ via $\inn_h$, so that $i_s$ is bijective, and
where $f_s$ is induced by the morphism
$(M_\QQ,\{s_0\})\to(\GL_{n,\QQ},Y)$ given by the inclusion of $M_\QQ$
in $\GL_{n,\QQ}$ followed by $\inn_h$, and
$s_0\mapsto\inn_h\circ\rho\circ s_0$.

The set $\GL_n(\QQ)\backslash(Y\times\GL_n(\AAf))/\GL_n(\Zhat)$ is the
set of isomorphism classes of $\ZZ$-Hodge structures $W$ such that
$W_\RR$ is isomorphic to~$V_{s_0,\RR}$ (to $(y,g)$ in
$Y\times\GL_n(\AAf)$ one associates the Hodge structure
$y\colon\SS\to\GL_{n,\RR}$ on the lattice~$\QQ^n\cap g\Zhat^n$). Its
subset of isomorphism classes of $W$ such that $W_\QQ$ is isomorphic
to~$V_{s_0,\QQ}$ is in bijection with
$S=C(\QQ)\backslash\GL_n(\AAf)/\GL_n(\Zhat)$, where $C_\QQ$ is the
centralizer in $\GL_{n,\QQ}$ of~$M_\QQ$. Being the centralizer of a
torus, $C_\QQ$ is connected and reductive (actually, $C_\Qbar$ is
isomorphic to a product of $\GL_{d_i,\Qbar}$'s). 
\vglue4pt
Note that $M(\AAf)$ acts on $S$ by left-multiplication. The image in
$S$ of $(\Zhat\otimes O_F)^*{\cdot}s$ is now simply the $(\Zhat\otimes
O_F)^*$-orbit of the class of $h$, where $(\Zhat\otimes O_F)^*$ acts
via $r_{\tilde{s}_0}\colon (\Zhat\otimes O_F)^*\to M(\AAf)$ and left
multiplication by~$M(\AAf)$. Let $L=\GL_n(\AAf)/\GL_n(\Zhat)$ be the
set of $\ZZ$-lattices in~$\QQ^n$ (or, equivalently, $\Zhat$-lattices
in~$\AAf^n$). The following lemma can be seen as a comparison between
the sizes of the $M(\Zhat)$-orbits $M(\Zhat)\ol{h}^S$ and
$M(\Zhat)\ol{h}^L$ of $h$ in $L$ and in~$S$.

\specialnumber{4.3.1}\proclaim{Lemma}\label{lem:SandL}
There exists an integer $m\geq1$ such that for all $h$ in
$\GL_n(\AAf)$ we have\/{\rm :}\/ 
$$
|M(\Zhat)\;\ol{h}^S| \geq |M(\Zhat)\;\ol{h}^L|, \quad\hbox{where
 $M(\Zhat)$ acts on $L$ via $m^{\rm th}$ powers.}
$$
\endproclaim

\demo{Proof}
The fact that $C(\RR)$ is not necessarily compact gives us some
trouble, and so we decompose $C_\QQ$ up to isogeny. There exists a
connected reductive subgroup $M'_\QQ$ of $C_\QQ$ such that we have a
short exact sequence:
\begin{equation} \speqnu{4.3.2}\label{eqn:isogeny}
1\lto F_\QQ \lto M_\QQ\times M'_\QQ \lto C_\QQ\lto 1,
\end{equation} 
with $F_\QQ=M_\QQ\cap M'_\QQ$ a finite group scheme. To get such an
$M'_\QQ$, consider the decomposition up to isogeny of $C_\QQ$ into its
center (that contains $M_\QQ$) and its semi-simple part, and use the fact that
the center decomposes up to isogeny into $M_\QQ$ and another
factor. All groups in (\ref{eqn:isogeny}) are closed subgroup schemes
of $\GL_{n,\QQ}$, which gives each of them a $\ZZ$-structure. There is
a nonempty open part $\Spec(\ZZ[1/N])$ of $\Spec(\ZZ)$ over which $F$
is finite \'etale, $M$ and $M'$ tori, $C$ reductive and the sequence
exact (for the \'etale topology). Let $e\geq 1$ be an integer that
annihilates $F_\QQ$ and $M(\QQ)\cap U$, where $U$ is the maximal
compact (open) subgroup of~$M(\AAf)$ (see
\cite[Prop.~3.16]{PlatonovRapinchuk}).

Let $h$ be in~$\GL_n(\AAf)$. Let $x$ be in $M(\Zhat)$ and suppose that
$x$ stabilizes $\ol{h}$ in~$S$. Then there exist $q$ in $C(\QQ)$ and
$k$ in $\GL_n(\Zhat)$ such that in $\GL_n(\AAf)$ we have: 
\begin{equation} \label{eqn:xh=qhk}
xh=qhk,\quad \hbox{i.e.,}\quad x=q{\cdot}hkh^{-1}.\speqnu{4.3.3}
\end{equation} 
Since $x$ and $q$ commute, we have: 
\begin{equation} \speqnu{4.3.4}\label{eqn:x^e}
x^e = q^e{\cdot}hk^eh^{-1}.
\end{equation} 
As $\rH^1(\Gal(\Qbar/\QQ),F(\Qbar))$ is annihilated by~$e$, there
exist $q_1$ in $M(\QQ)$ and $q_2$ in $M'(\QQ)$ such that $q^e=q_1q_2$
in~$C(\QQ)$. By (\ref{eqn:xh=qhk}), $hkh^{-1}$ is in~$C(\AAf)$. Since
the sequence (\ref{eqn:isogeny}) extends as indicated over~$\ZZ[1/N]$,
there are $k_1$ and $k_2$ in $\GL_n(\AAf)$ such that: 
\begin{equation} \speqnu{4.3.5}
hk^eh^{-1} = hk_1h^{-1}{\cdot}hk_2h^{-1},
\end{equation} 
with $hk_1h^{-1}$
in $M(\AAf)$ and $hk_2h^{-1}$ in $M'(\AAf)$. 
Rewriting (\ref{eqn:x^e}) gives: 
\begin{equation} \speqnu{4.3.6}
x^e = q_1q_2{\cdot}hk_1h^{-1}{\cdot}hk_2h^{-1} = 
q_1{\cdot}hk_1h^{-1}{\cdot}q_2hk_2h^{-1}.
\end{equation} 
It follows that $q_2hk_2h^{-1}$ is in $F(\AAf)$, hence that:
\begin{equation} \speqnu{4.3.7}
x^{e^2} = q_1^e{\cdot}hk_1^eh^{-1}, \quad\hbox{in $M(\AAf)$.}
\end{equation} 
This identity shows that $q_1^e$ is in $U\cap M(\QQ)$. Hence $q_1$ is
in $U\cap M(\QQ)$, and since $e$ annihilates $U\cap M(\QQ)$, we have: 
\begin{equation} \speqnu{4.3.8}
x^{e^2} = hk_1^eh^{-1}, \quad\hbox{in $M(\AAf)$.}
\end{equation} 
We conclude that for $m=e^2$, we have, for $x$ in $M(\Zhat)$
stabilizing $\ol{h}^S$, that $x^m$ stabilizes $\ol{h}^L$ (for the
usual action of $M(\AAf)$ on~$L$). 
\enddemo

Lemma~\ref{lem:SandL} implies that in order to prove
Theorem~\ref{thm3.1} we may as well prove the lower bound in question
for the action of $M(\Zhat)$ on the set $L$, if we replace the
reciprocity morphism~$r_{\tilde{s}_0}$ by $mr_{\tilde{s}_0}$, for a
suitable~$m\geq1$. Note that $mr_{\tilde{s}_0}$ is still surjective,
as a morphism of tori over~$\QQ$. 

For $h$ in $\GL_n(\AAf)$, its class in $L$ is a $\ZZ$-lattice $V_h$
in~$\QQ^n$, and $\QQ^n$ is equipped with an action of~$M_\QQ$. Let
$\MT(V_h)$ be the scheme theoretic closure of $M_\QQ$
in~$\bGL(V_h)$. In this situation, Lemma~\ref{lem:condatp} says that
$\MT(V_h)_{\FF_p}$ is a torus if and only if $M$ fixes~$V_h$.

We have reduced the proof of Theorem~\ref{thm3.1} to the following
statement, to be proved in the next section.

\specialnumber{4.3.9}\proclaim{Proposition}\label{prop:T-orbit}
Let $n$ be a positive integer.  Let $T$ be a torus over $\ZZ[1/n]${\rm ,}
acting on a free $\ZZ[1/n]$\/{\rm -}\/module $V$ of finite rank. Let
$$S=\GL(V_{\AAf})/\GL(V_{\Zhat[1/n]})$$ be the set of
$\ZZ[1/n]$-lattices in~$V_\QQ$. Let $T(\Zhat[1/n])$ act on $S$ via
left multiplications. For $W$ in $S${\rm ,} let $P_W$ be the set of primes
$p$ that do not divide $n$ and such that $T_{\ZZ_p}$ does not
fix~$W_{\ZZ_p}$. Then there exists a positive real number $c$ such
that for each $W$ in $S$ we have\/{\rm :}
$$
|T(\Zhat[1/n]){\cdot}W| \geq \prod_{p\in P_W}cp.
$$
Equivalently{\rm ,} as $S$ is a restricted product over the primes not
dividing~$n${\rm ,} there is a positive real number $c$ such that for all
$p$ not dividing~$n${\rm ,} and all $\ZZ_p$-lattices $W$ in $V_{\QQ_p}$ that
are not fixed by~$T_{\ZZ_p}${\rm ,} we have\/{\rm :}\/ 
$$
|T(\ZZ_p){\cdot}W| \geq cp.
$$
\endproclaim

4.4. {\it Proof of Proposition}~\ref{prop:T-orbit}.
  Let $K$ be a
splitting field of $T_\QQ$, and let $X$ be the group of characters
of~$T_K$. Let $p$ be a prime number that does not divide~$n$, and let
$W$ be a $\ZZ_p$-lattice in $V_{\QQ_p}$ that is not fixed
by~$T_{\ZZ_p}$. Let $W'$ be the largest sublattice of $W$ that is
fixed by $T_{\ZZ_p}$, as constructed in the last part of the proof of
Lemma~\ref{lem:condatp}. Let $L$ be the kernel of multiplication by
$p$ on~$W/W'$. We view $L$ as a sub $\FF_p$-vector space of
$W'_{\FF_p}$ via $L=(W\cap p^{-1}W')/W'=(pW\cap W')/pW'$, where the
second equality comes from multiplication by~$p$. 

An element of $T(\ZZ_p)$ that stabilizes~$W$  also
stabilizes~$L$. Hence $T(\ZZ_p){\cdot}W$ has at least as many elements
as~$T(\FF_p){\cdot}L$, where we let $T(\FF_p)$ act on the set of
subspaces of~$W'_{\FF_p}$. By construction, $L$ is nonzero, and has
intersection zero with $W'_{\Fbar_p,x}$ for each $x$ in~$X$, where
$W'_{\Fbar_p}=\oplus_x W'_{\Fbar_p,x}$ is the $X$-grading
corresponding to the $T_{\FF_p}$-action on~$W'_{\FF_p}$. This implies
that~$L$, viewed as an $\FF_p$-valued point in some Grassmannian, is
not fixed by the action of~$T_{\FF_p}$. 

Let $T'_{\FF_p}$ be the stabilizer of $L$
in~$T_{\FF_p}$. Lemma~\ref{lem:concomp} below says that the order of
the group of connected components of $T'_{\Fbar_p}$ is bounded
independently of $p$ and~$W$. Put
$T''_{\FF_p}:=T_{\FF_p}/T'_{\FF_p}$. Then $T''_{\FF_p}$ is a
nontrivial torus over~$\FF_p$, and by Lemma~\ref{lem:exseqtori}, the
morphism $T(\FF_p)\to T''_{\FF_p}(\FF_p)$ has its cokernel of order
bounded independently of $p$ and~$W$. The proof of
Proposition~\ref{prop:T-orbit} is now finished   if we note that
$T''_{\FF_p}(\FF_p)$ has at least $p{-}1$ elements (see
\cite[Lemma~3.5]{Nori1}).\hfill\qed

\specialnumber{4.4.1}\proclaim{Lemma}\label{lem:concomp}
Let $k$ be an algebraically closed field{\rm ,} and $T$ a $k$\/{\rm -}\/torus. Let $V$
be a finite-dimensional $k$\/{\rm -}\/vector space with an action by~$T$. Then
the set of stabilizers $T_W${\rm ,} for $W$ running through the set of
subspaces of $V${\rm ,} is finite. The set of groups of connected components
of these stabilizers{\rm ,} up to isomorphism{\rm ,} is finite{\rm ,} and bounded in
terms of the dimension of\/ $V$ and the set of characters of $T$ that
do occur in~$V$.
\endproclaim

\demo{Proof}
Let us consider the set $S$ of subspaces $W$ of a fixed dimension,
call it~$d$. Then we have a natural injection from $S$ into
$\PP(\Lambda^d(V))$, compatible with the $T(k)$-action. The image of a
$W$ under this map is the line generated by $w_1\wedge\cdots\wedge
w_d$, where $w$ is any $k$-basis of~$W$. Hence the set of stabilizers
of the elements of $S$ is contained in the set of stabilizers of
elements of $\PP(\Lambda^d(V))$. This reduces the proof of the lemma
to the case of one-dimensional subspaces (we replace $V$
by~$\Lambda^d(V)$).

Let $X$ be the character group of $T$, and let $V=\oplus_x V_x$
be the $X$-grading of $V$ given by the $T$-action. Of course, almost
all $V_x$ are zero. For $v$ in $V$, we have $v=\sum_x v_x$,
and we let $\Supp(v)$ be the set of $x$ with $v_x\neq0$. For $v$
nonzero, the stabilizer in $T$ of the element $kv$ in $\PP(V)$ is the
intersection of the kernels of the $x-x'$ with $x$ and
$x'$ in~$\Supp(v)$. Since the set of such differences is finite,
the claim follows. 
\enddemo

\specialnumber{4.4.2}\proclaim{Lemma}\label{lem:exseqtori}
Let $T'$ be the kernel of a surjective morphism $T\to T''$ of tori
over~$k${\rm ,} with $k$ a finite field. Let $\Phi$ be the group of
connected components of~$T'_{\ol{k}}$. Then there are  the exact
sequence\/{\rm :}
$$
0\lto T'(k)\lto T(k) \lto T''(k)\lto \rH^1(\Gal(\ol{k}/k),\Phi), 
$$
and the upper bound\/{\rm :} 
$$
|\rH^1(\Gal(\ol{k}/k),\Phi)| \leq |\Phi|.
$$
\endproclaim

\demo{Proof}
One has, of course, the long exact sequence coming from Galois
cohomology:
$$
0\lto T'(k)\lto T(k) \lto T''(k)\lto \rH^1(\Gal(\ol{k}/k),T'(\ol{k})). 
$$
We combine this with the sequence:
$$
\rH^1(\Gal(\ol{k}/k),(T')^0(\ol{k})) \lto
\rH^1(\Gal(\ol{k}/k),T'(\ol{k})) \lto \rH^1(\Gal(\ol{k}/k),\Phi)
$$
coming from the short exact sequence: 
$$
0\lto (T')^0\lto T'\lto T'/(T')^0\lto 0. 
$$
Lang's Theorem ([10,  Thm.~6.1]) implies that
$\rH^1(\Gal(\ol{k}/k),(T')^0(\ol{k}))=0$. The upper bound for
$|\rH^1(\Gal(\ol{k}/k),\Phi)|$ follows from the fact that this
cohomology group is just the group of coinvariants for the action of
$\Gal(\ol{k}/k)$ on~$\Phi$.
\enddemo

\section{Images under Hecke correspondences}\label{sec4}

In this section we prove that the images under irreducible components
of certain Hecke correspondences of an irreducible Hodge generic
subvariety of a Shimura variety defined by a semi-simple algebraic
group of adjoint type are irreducible.

\proclaim{Theorem}\label{thm4.1}
Consider a Shimura variety defined by a Shimura datum $(G,X)$ where
$G$ is a semi\/{\rm -}\/simple algebraic group of adjoint type. Let $K$ be a neat
compact open subgroup of\/ $G(\AAf)$ that is the product of compact
open subgroups $K_p$ of~$G(\QQ_p)$. Let $X^+$ be a connected component
of $X$ and let $S$ be the image of $X^+\ti\{1\}$ in
$\Sh_K(G,X)_\CC$. Let $Z$ be an irreducible Hodge generic subvariety
of $S$ containing a nonsingular special point. Then there exists a
nonzero integer $n$ such that if $q$ is an element of
$G(\QQ)^+:=G(\QQ)\cap G(\RR)^+$ whose image in $G(\QQ_l)$ is in $K_l$
for every $l$ dividing $n$, then $T_q(Z)$ is irreducible{\rm ,} with $T_q$
the correspondence on $S$ given by the action of $q$ on~$X^+$.
\endproclaim 

\demo{Proof}
By \cite[2.1.2]{De2} we have $S=\Gamma\backslash X^+$ with $\Gamma =
G(\QQ)^+ \cap K$.  For $q$ in~$G(\QQ)^+$ the correspondence $T_q$ is
defined as follows. Consider the diagram:
$$
S\stackrel{\pi}{\lfrom} X^+ \stackrel{q\cdot}{\lto} X^+
\stackrel{\pi}{\lto} S,
$$
with $\pi$ the quotient map for the action of~$\Gamma$. The morphism
$\pi\circ q\cdot$ is the quotient for the action of $q^{-1}\Gamma q$;
hence $\pi$ and $\pi\circ q\cdot$ both factor through the quotient
$S_q=\Gamma_q\backslash X^+$ with $\Gamma_q=\Gamma\cap q^{-1}\Gamma
q$. With this notation, $T_q$ is the correspondence: 
$$
S\stackrel{\pi_{q,1}}{\lfrom} S_q\stackrel{\pi_{q,2}}{\lto} S.
$$
In particular, $T_qZ=\pi_{q,2}(\pi_{q,1}^{-1}Z)$. In order to show
that $T_qZ$ is irreducible, it suffices to show that
$Z_q:=\pi_{q,1}^{-1}Z$ is. Since $\pi_{q,1}$ is a covering,
restriction to the smooth loci in $Z$ and $Z_q$ gives a covering
$\pi_{q,1}\colon Z_q^\sm\to Z^\sm$. Since $Z_q^\sm$ is
Zariski-dense in $Z_q$, it is enough to prove that $Z_q^\sm$ is
irreducible.

In order to find an integer $n$ as in the theorem, we choose a
faithful representation $\xi$ of $G$ on a $\QQ$-vector space~$V$. Then
$\xi$ gives a polarizable variation of Hodge structures on the
constant sheaf $V_X$ (\cite[1.8]{Moonen1}). Since $K$ is neat,
$\Gamma$ acts freely on $X^+$, hence $V_{X^+}$ descends to~$S$. Let
$V_{\Zhat}$ be a $K$-invariant $\Zhat$-lattice in~$V_{\AAf}$, and let
$V_\ZZ=V_{\Zhat}\cap V$. Then $V_\ZZ$ is a $\Gamma$-invariant lattice
in~$V$. This gives us a polarizable variation $\calV_\ZZ(\xi)$ of
$\ZZ$-Hodge structure on $S$, hence on~$Z^\sm$. Let $s$ be a Hodge
generic point in $Z^\sm$ and $x$ a point of $X^+$ lying above it. This
gives an isomorphism between the fiber $\calV_\ZZ(\xi)_s$ of
$\calV_\ZZ(\xi)$ at $s$ and $V_\ZZ$, such that the Mumford-Tate group
of $\calV_\ZZ(\xi)_s$ corresponds to~$\xi(G)$.

On the other hand, the fundamental group $\pi_1(Z^\sm,s)$ acts on
$\calV_\ZZ(\xi)_s$ and hence on $V_\ZZ$; let $\Pi$ be its image
in~$\GL(V_\ZZ)$. Then $\Pi$ is a finitely generated subgroup
of~$\GL(V_\ZZ)$. By a theorem of Andr\'e (see [1, Th.~1.4]), the connected component of the Zariski closure
of
$\Pi$ in the algebraic group $\GL(V)$ is~$\xi(G)$. (Here we use the
assumption that $Z$ contains a nonsingular special point.) Since
$\Gamma$ is an arithmetic subgroup of $G(\QQ)$, it is Zariski dense
in~$G$, and the Zariski closure of $\xi(\Gamma)$ is~$\xi(G)$.

Since $\calV_\ZZ(\xi)$ exists over $S$, the image $\Pi$ of
$\pi_1(Z^\sm,s)$ is contained in that of $\pi_1(S,s)$, i.e.,
in~$\xi(\Gamma)$. The inclusion
$\Pi\subset\xi(\Gamma)\stackrel{\xi^{-1}}{\to}\Gamma$ makes $\Gamma$
into a $\pi_1(Z^\sm,s)$-set. For $q$ in $G(\QQ)^+$, the
$\pi_1(Z^\sm,s)$-set that corresponds to the covering $\pi_{q,1}\colon
Z_q^\sm\to Z^\sm$ is isomorphic to $\Gamma /q^{-1}\Gamma q\cap\Gamma$;
hence $Z^\sm$ is connected (and hence irreducible) if and only if
$\Gamma /q^{-1}\Gamma q\cap\Gamma$ is transitive. Now since
$\pi_1(Z^\sm,s)$ and $\Gamma$ have the same Zariski closure in the
$\ZZ$-group scheme $\GL(V_\ZZ)$, Nori's Theorem~\ref{thmNori} below
gives a nonzero integer $n$, such that for all nonzero integers $m$
prime to $n$, $\pi_1(Z^\sm,s)$ and $\Gamma$ have the same image
in~$\GL(V_{\ZZ/m\ZZ})$. (Indeed, $G(\CC)$ is a semisimple complex Lie
group, hence $\pi_1(G(\CC))$ is finite by~\cite[Thm.~2(c)]{Serre1}.)
Now let $q$ be in $G(\QQ)^+$, such that $q_l$ is in $K_l$ for all $l$
dividing~$n$. Then all we have to show is that there is a nonzero
integer $m=\prod_ll^{m_l}$, prime to $n$, such that $\Gamma\cap
q^{-1}\Gamma q$ contains the kernel of the natural map $\Gamma \to
\GL(V_{\ZZ/m\ZZ})$. We have:
$$
q^{-1}\Gamma q\cap\Gamma = q^{-1}K q\cap K\cap G(\QQ)^+ = 
\left(\prod_l q^{-1}K_l q\cap K_l\right)\cap G(\QQ)^+.
$$
For $l$ with $q_l\in K_l$ (in particular, for $l$ dividing $n$), we
have $q^{-1}K_l q\cap K_l=K_l$, and we put $m_l=0$. For the remaining
finitely many $l$, we take $m_l$ sufficiently large such that:
$$
\ker\left(K_l\lto\GL(V_{\ZZ/l^{m_l}\ZZ})\right) \subset q^{-1}K_lq.
$$
Then we have the last inclusion for all $l$, so by taking the product
over all $l$ we have: 
$$
\ker\left(K\lto\GL(V_{\ZZ/m\ZZ})\right) \subset q^{-1}Kq,
\quad\hbox{hence}\quad 
\ker\left(\Gamma \lto\GL(V_{\ZZ/m\ZZ})\right)\subset q^{-1}\Gamma q.
$$
\phantom{wind}
\enddemo
 
\proclaimtitle{Nori (Thm.~5.3 of \cite{Nori1})}
\proclaim{Theorem}\label{thmNori}
Let $H$ be a finitely generated subgroup of $\GL_n(\ZZ)$. Let $\ol{H}$
be its Zariski closure in $\GL_{n,\ZZ}$. Suppose that $\ol{H}(\CC)$
has a finite fundamental group{\rm ;} then the closure of $H$ in
$\GL_n(\Zhat)$ is open in the closure of $\ol{H}(\ZZ)$
in~$\GL_n(\Zhat)$.
\endproclaim 

\numbereddemo{{R}emark}
A careful reader might be worried by our use of Theorem~5.3 of Nori's
article \cite{Nori1}, because Theorem~5.2 of the same article is
clearly wrong. (The image of $\SL_2(\ZZ)$ in $\GL_3(\QQ)$ via the
symmetric square of the standard representation gives a counterexample.) Probably, the problem with Theorem~5.2 is
only typographical: should $[A,A](\hat{R})$ be replaced by
$[A,A](R)^\wedge$, the closure of $[A,A](R)$ in $\GL_n(\hat{R})$? Or
by $[A(R),A(R)]^\wedge$?  Anyway, Theorem~5.3 can be deduced from
Theorem~5.4, which in turn can be found in other references
(\cite[Thm.~7.14]{PlatonovRapinchuk},~\cite{Pink1}).
\enddemo

 \vglue-12pt
\section{Density of Hecke orbits}\label{sec5}
\vglue-8pt

We will now prove the following result about the density of Hecke
orbits in Shimura varieties.

\proclaim{Theorem}\label{thm5.1}
Let $(G,X)$ be a Shimura datum with $G=G_1\times\cdots\times G_r$ a
semi\/{\rm -}\/simple algebraic group of adjoint type with simple
factors~$G_1,\ldots,G_r$.  Let $K$ be a compact open subgroup of\/
$G(\AAf)$ that is the product of compact open subgroups $K_p$
of~$G(\QQ_p)$.  Let $X^+$ be a connected component of~$X${\rm ,} let $S$ be
the image of $X^+ \ti\{1\}$ in \/ $\Sh_K(G,X)_\CC${\rm ,} and put $\Gamma =
G(\QQ)^+\cap K$. Then there exists an integer $n$ such that for all
prime numbers $p\geq n$ the following holds.

Let $q$ be an element of $G(\QQ_p)$ such that for all $j$ the
projection of $q$ in $G_j(\QQ_p)$ is not contained in a compact
subgroup. The connected components of the correspondence that $T_q$
induces on $S$ are the~$T_{q_i}$, induced by $q_i$ in $G(\QQ)^+$
acting on~$X^+${\rm ,} such that\/{\rm :}\/
$$
G(\QQ)^+ \cap K q K = \coprod_i \Gamma q_i^{-1}\Gamma.
$$
Then{\rm ,} for all $i$ and for all $s$ in~$S${\rm ,} the
$T_{q_i}+T_{q_i^{-1}}$-orbit $\cup_{n\geq0}(T_{q_i}+T_{q_i^{-1}})^ns$
is dense in $S$ for the Archimedean topology.
\endproclaim 

\demo{Proof}
Let $p_j\colon G\to G_j$ denote the $j^{\rm th}$ projection. For each~$j$,
let $n_j$ be as in Proposition~\ref{prop5.2} applied to $G_j$
and~$p_j(K)$, and let $n$ be the maximum of the~$n_j$. Let $p\geq n$
be prime, and let $q$ in $G(\QQ_p)$ and the $q_i$ in $G(\QQ)$ be as in
the theorem. For each~$i$, let $\Gamma_i$ denote the subgroup of
$G(\QQ)^+$ generated by $\Gamma$ and~$q_i$. Then we have
$S=\Gamma\backslash X^+$, and for $s$ in $S$, the
$T_{q_i}+T_{q_i^{-1}}$-orbit of $s$ is the image in $S$ of a
$\Gamma_i$-orbit in~$X^+$. Hence it suffices to prove that each
$\Gamma_i$ is dense in~$G(\RR)^+$.

Let $H$ be the closure in $G(\RR)^+$ of say ~$\Gamma_1$. Then $H$ is a
Lie subgroup of~$G(\RR)$. We let $H^+$ be the connected component of
the identity. Now $\Gamma$ normalizes $H^+$, hence its Lie
algebra. Since $\Gamma$ is Zariski dense in $G_\RR$
(\cite[Thm.~4.10]{PlatonovRapinchuk}; note that no $G_i(\RR)$ is
compact), this implies that $H^+$ is a product of simple factors
of~$G(\RR)^+$. We claim that if $H^+$ contains a simple factor of
$G_i(\RR)^+$, then it contains all of~$G_i(\RR)^+$. This follows from
two facts. One: $H^+\cap G(\QQ)$ is dense in~$H^+$. Two: for every
$i$, and for every factor $G_{i,j}$ of $G_{i,\Qbar}$ the map
$G_i(\QQ)\to G_{i,j}(\Qbar)$ is injective. The second fact is a direct
consequence of the statement (\cite[6.21(ii)]{BorelTits1}) that
$G_i=\Res_{F/\QQ}G_i'$ for some finite extension $F$ of $\QQ$ and some
absolutely simple group $G_i'$ over~$F$.

Let now $I$ be the subset of $\{1,\ldots,r\}$ of the $i$ such that
$H^+$ does not contain $G_i(\RR)^+$, let $G_I$ be the product of the
$G_i$ with $i$ in $I$, and let $f\colon G(\RR)^+\to G_I(\RR)^+$ be the
projection. Then $H^+$ is the kernel of $f$, hence $f(H)$ is discrete
in~$G_I(\RR)$. But then $f(\Gamma_1)$, being a subgroup of $f(H)$, is
discrete. Now suppose that $I$ is not empty. By
Proposition~\ref{prop5.2}, $f(\Gamma)$ has infinite index in
$f(\Gamma_1)$, which contradicts the fact that $f(\Gamma)\backslash
G_I(\RR)^+$ has finite volume (see
\cite[Thm.~4.13]{PlatonovRapinchuk}). Hence $I$ is empty, and
$H=G(\RR)^+$.
\enddemo 

\proclaim{Proposition}\label{prop5.2}
Let $G$ be a simple algebraic group {\rm (}\/of adjoint type\/{\rm )} over~$\QQ$ with
$G(\RR)$ not compact. Let $K$ be a compact open subgroup of $G(\AAf)$
and put $\Gamma=G(\QQ)^+\cap K$. Then there is an integer $n$ such
that for all prime numbers $p\geq n$ the following holds.

When $g$ is an element of $G(\QQ_p)$ that is not contained in a compact
subgroup{\rm ,} then for all $q$ in $G(\QQ)\cap KgK${\rm ,} $\Gamma$ has infinite
index in the group generated by $\Gamma$ and~$q$.
\endproclaim

{\it Proof}.
Let $V$ be a finite-dimensional absolutely irreducible faithful
representation of~$G_E$, where $E$ is a suitable finite extension
of~$\QQ$. As $\Gamma$ is Zariski dense in $G_\Qbar$
(\cite[Thm.~4.10]{PlatonovRapinchuk}), $V$ is an absolutely
irreducible representation of~$\Gamma$. By the double centralizer
theorem (\cite[XVII, Cor.~3.5]{Lang2}), $\End_E(V)$ is spanned
(over~$E$) by the images of the $\gamma$ in~$\Gamma$. Let $O_E$ be the
ring of integers of~$E$, and let $V_{O_E}$ be a $\Gamma$-invariant
$O_E$-lattice in~$V$. Let $\gamma_1,\ldots,\gamma_m$ be in $\Gamma$
whose images $\ol{\gamma_1},\ldots,\ol{\gamma_m}$ span
$\End_E(V)$. Let $n$ be the index of
$O_E\ol{\gamma_1}+\cdots+O_E\ol{\gamma_m}$ in
$\End_{O_E}(V_{O_E})$. For $p$ prime put $O_{E,p}=\ZZ_p\otimes O_E$,
and $E_p=\QQ_p\otimes E$. For $p>n$ we have
$O_{E,p}\ol{\gamma_1}+\cdots+O_{E,p}\ol{\gamma_m}=
\End_{O_{E,p}}(V_{O_{E,p}})$.

Let now $p>n$ be prime, and let $g$ be an element of $G(\QQ_p)$ that
is not contained in a compact group. Then not all eigenvalues of $g$
acting on $\Qbar_p\otimes_\QQ V$ (as $\Qbar_p$-vector space) have
absolute value~$1$. Replacing $g$, $q$ and $q'$ by their inverses if
necessary, we may suppose that at least one of the eigenvalues of $g$
has absolute value~$>1$. Then the matrix of $g$ with respect to an
$O_{E,p}$-basis of $V_{O_{E,p}}$ has at least one coefficient that is
not in~$O_{E,p}$.

Let $q$ be an element of~$KgK$. Then we have, in~$G(\QQ_p)$,
$q=k_1gk_2$, with $k_1$ and $k_2$ in the image $K_p$ of~$K$. So in
$\GL_{E_p}(V_{E_p})$ we have $\ol{q}=\ol{k_1}\ol{g}\ol{k_2}$ with
$\ol{k_1}$ and $\ol{k_2}$ in~$\GL_{O_{E,p}}(V_{O_{E,p}})$. But then we
have $\ol{g}=\ol{k_1}^{-1}\ol{q}\ol{k_2}^{-1}$. It follows that all
$u\ol{g}v$ with $u$ and $v$ in $\End_{O_{E,p}}(V_{O_{E,p}})$ are
$O_{E_p}$-linear combinations of elements of the form
$\ol{\gamma}\ol{q}\ol{\gamma'}$ with $\gamma$ and $\gamma'$
in~$\Gamma$. Hence there are  $\gamma$ and $\gamma'$ such that the
trace of $q'=\gamma q\gamma'$ is not in~$O_{E,p}$. But then $q'$ has
an eigenvalue (in~$\Qbar_p$) of absolute value $>1$; hence $q'$ is not
contained in a compact subgroup of~$G(\QQ_p)$. If $q$ is in
$G(\QQ)\cap KgK$, then we get a $q'$ in $\Gamma q\Gamma$ where no
nontrivial power lies in~$\Gamma$.
\hfill\qed
\vglue6pt

One would like to have a generalization of Proposition~\ref{prop5.2}
that requires only a weaker hypothesis on $q$ and that gives a better
understanding of the primes $p$ that are to be excluded. The following
proposition, that will not be used in the rest of this article, gives
such a result in the simply connected case. It would be useful to have
a version of it in the adjoint case (maybe the results
in~\cite{Margulis1} can be used here).

\proclaim{Proposition}\label{prop5.3}
Let $G$ be a semisimple algebraic group over $\QQ$ whose adjoint is
simple. Suppose that $G(\RR)$ is not compact. Let $K$ be a compact
open subgroup of~$G(\AAf)$. For $p$ prime let $K_p$ be the image of
$K$ under projection to~$G(\QQ_p)$. Let $p$ be a prime for which $K_p$
is a maximal compact subgroup of~$G(\QQ_p)$ and let $q$ be an element
of $G(\QQ)$ such that the image of $q$ in $G(\QQ_p)$ is not
in~$K_p$. Then $\Gamma=G(\QQ)\cap K$ is of infinite index in the group
$\Gamma_q$ generated by $\Gamma$ and~$q$.
\endproclaim

{\it Proof}.
By strong approximation (\cite[Thm.~7.12]{PlatonovRapinchuk}),
$G(\QQ)$ is dense in~$G(\AAf)$. Hence $\Gamma$ is dense
in~$K$. Suppose that $\Gamma_q/\Gamma$ is finite. Let $\ol{\Gamma}$
and $\ol{\Gamma_q}$ be the closures of $\Gamma$ and $\Gamma_q$
in~$G(\QQ_p)$. We have $\ol{\Gamma}=K_p$. Since $\Gamma_q$ is a finite
union of $\Gamma$-cosets, $\ol{\Gamma_q}$ is a finite union of
$\ol{\Gamma}$-cosets. Hence $\ol{\Gamma_q}$ is compact. As $q$ is not
in~$K_p$, and~$K_p$ is a maximal compact open subgroup of~$G(\QQ_p)$,
we have a contradiction.\hfill\qed

 \vglue-8pt
\section{Proof of the main result}\label{sec6}
 \vglue-9pt

The aim of this section is to prove Theorem~1.2. So let
$(G,X)$ be a Shimura datum, let $K$ be a compact open subgroup of
$G(\AAf)$, let $V$ be a faithful finite-dimensional representation
of~$G$ and let $Z$ be a closed irreducible curve in $\Sh_K(G,X)_\CC$
that contains an infinite set $\Sigma$ of special points $x$ such that
all $[V_x]$ with $x$ in $\Sigma$ are equal. 

As we have seen in Section~\ref{sec2}, it suffices to prove the result
in the case where $G$ is semi-simple and adjoint, $K$ is neat and $Z$
Hodge generic.  Given these conditions, we write $G=G_1\times\cdots\times
G_r$ with the $G_i$ simple. We will also suppose that $K$ is the
product of compact open subgroups $K_p$ of the~$G(\QQ_p)$. We let $S$
be the connected component of $\Sh_K(G,X)$ that contains~$Z$.  We
assume that there is a $\ZZ$-structure on~$V$ given by the choice of a
$K$-invariant $\Zhat$-lattice in~$V_{\AAf}$. Then $V$ induces a
variation of $\ZZ$-Hodge structure on~$S$. For each $s$ in $S$ we have
its Mumford-Tate group $\MT(V_s)$ which is a closed subgroup scheme
of the $\ZZ$-group scheme~$\bGL(V_s)$.

Let $X^+$ be a connected component of~$X$. After replacing $Z$ by an
irreducible component of its image under a suitable Hecke
correspondence, we may suppose that $S$ is the image of
$X^+\times\{1\}$ in~$\Sh_K(G,X)$.

\proclaim{Theorem}\label{thm6.1}
Assume that $p$ is a prime and $m$ an element of $G(\QQ_p)${\rm ,} such
that\/{\rm :}
\begin{itemize}
\ritem{1.} $p$ does not divide the integer $n$ of Theorem~{\rm \ref{thm4.1},}
applied to the subvariety $Z$ of\/ $\Sh_K(G,X)${\rm ;}
\ritem{2.} $p\geq n$ with $n$ as in Theorem~{\rm \ref{thm5.1}} applied to $(G,X)${\rm ,}
$K$ and $X^+${\rm ,} and for every $\QQ$-simple factor $G_i$ of $G${\rm ,} the image
of $m$ in $G_i(\QQ_p)$ is not contained in a compact subgroup\/{\rm ;}\/
\ritem{3.} $Z\subset T_mZ$, with $T_m$ the correspondence induced by $m$
on $\Sh_K(G,X)$.
\end{itemize}
Then $Z=S${\rm ;} hence $Z$ is of Hodge type.
\endproclaim 

\demo{Proof}
Assume that $p$ and $m$ satisfy the three conditions. We have
$S=\break\Gamma\backslash X^+$, with $\Gamma$ the intersection of $K$
and~$G(\QQ)^+$. The connected components of the correspondence that
$T_m$ induces on $S$ are the $T_{m_i}$, induced by $m_i$ in
$G(\QQ)^+$ acting on $X^+$, such that:
$$
G(\QQ)^+ \cap K m K = \coprod_i \Gamma m_i^{-1}\Gamma.
$$
Let $q$ be one of the $m_i$ such that $Z\subset T_q Z$. Since $q^{-1}$
is in $KmK$, and $m$ in $G(\QQ_p)$, the image of $q$ in $G(\QQ_l)$ is
in $K_l$ for all $l\neq p$, in particular, for all $l$
dividing~$n$. By Theorem~\ref{thm4.1}, $T_qZ$ and $T_{q^{-1}}Z$ are
irreducible. Hence $Z=T_qZ$, and $T_{q^{-1}}Z=T_{q^{-1}}T_qZ\supset
Z$, so that $T_{q^{-1}}Z=Z$. By Theorem~\ref{thm5.1}, all
$T_q+T_{q^{-1}}$-orbits in $S$ are dense; hence $Z=S$.
\enddemo

For each $s$ in $\Sigma$ we choose an element $\tilde{s}$ of $X^+$
such that $s=\ol{(\tilde{s},1)}$. Then   for each $s$ in
$\Sigma$ the Mumford-Tate group $\MT(V_s)=\MT(\tilde{s})$ is a closed
subgroup scheme of~$\bGL(V)$. We give $G$ a $\ZZ$-structure $G_\ZZ$
by taking its Zariski closure in~$\bGL(V)$.

Since all $V_{s,\QQ}$ for $s$ in $\Sigma$ are isomorphic, all
$\MT(\tilde{s})_\QQ$ are isomorphic. Let $F\subset\CC$ be a finite
extension of $\QQ$ that splits all these tori. Then $F$ contains the
reflex field of~$(G,X)$.  By taking $F$ large enough, we may and do
assume that $Z$ is defined over $F$ (as an absolutely irreducible
closed $F$-subscheme $Z_F$ of $\Sh_{K}(G,X)_F$).

In the rest of the proof of Theorem~1.2 we will distinguish
two cases, depending on the behavior of the $\MT(V_s)$ for $s$
in~$\Sigma$.
 
\numbereddemo{Definition}\label{def6.2}
For $s$ in $\Sigma$, let $i(s)$ be the number of prime numbers $p$
such that $\MT(V_s)_{\FF_p}$ is not a torus.
\enddemo

This function $i\colon\Sigma\to\ZZ$ depends on the choice of the
$\ZZ$-structure on~$V$, but another choice would give a function that
differs from $i$ by a bounded function. Also, replacing $Z$ by an
irreducible component of the image of $Z$ under a Hecke correspondence
only changes the function $i$ by a bounded function. The two cases
that we will distinguish depend on whether or not $i$ is bounded.

\demo{{\rm 7.3.} The case where $i$ is bounded}
In this section we assume that the function $i:\Sigma\to\ZZ$ defined
above is bounded. It follows that for all but finitely many prime
numbers~$p$, the subset $\Sigma_p$ of $s$ in $\Sigma$ with
$\MT(V_s)_{\FF_p}$ a torus is Zariski dense in~$Z$. Indeed, if $p$ is
such that $\Sigma_p$ is not dense, then we can replace $\Sigma$ with
the complement of its subset~$\Sigma_p$. If $i$ is bounded by~$B$,
then, for each $s$ in~$\Sigma$, there are at most $B-1$ primes other
than $p$ for which $\MT(V_s)_{\FF_p}$ is not a torus, etc. 
\enddemo

\specialnumber{7.3.1}\proclaim{Proposition}
Let $p$ be a prime such that $G_{\FF_p}$ is smooth over~$\FF_p$. Then
the set of subtori $\MT(\tilde{s})_{\ZZ_p}$ of $G_{\ZZ_p}$ for $s$ in
$\Sigma_p$ meets only a finite number of $G_\ZZ(\ZZ_p)$-conjugacy
classes of subtori of~$G_{\ZZ_p}$.
\endproclaim

\demo{Proof}
Let $p$ be prime. Let $s_1$ and $s_2$ be in~$\Sigma_p$. Let $h$ in
$\GL(V_\QQ)$ be an isomorphism from $(V_\QQ,\tilde{s}_1)$
to~$(V_\QQ,\tilde{s}_2)$. Then $h$ induces an isomorphism of
$\ZZ$-Hodge structures from $(h^{-1}V,\tilde{s}_1)$
to~$(V,\tilde{s}_2)$. Hence $h^{-1}V_{\ZZ_p}$ is a $\ZZ_p$-lattice in
$V_{\QQ_p}$ such that the Zariski closure of
$\MT(\tilde{s}_1)_{\QQ_p}$ in $\bGL(h^{-1}V_{\ZZ_p})$ is a
torus. Lemma~\ref{lem:condatp} (applied to
$T=\MT(V,\tilde{s}_1)_{\ZZ_p}$ and the free $\ZZ_p$-module
$h^{-1}V_{\ZZ_p}$) gives an element $c$ in $\GL(V_{\QQ_p})$
centralizing $\MT(V,\tilde{s}_1)_{\QQ_p}$ such that
$h^{-1}V_{\ZZ_p}=cV_{\ZZ_p}$. Hence there exists $k$ in
$\GL(V_{\ZZ_p})$ such that $h^{-1}=ck$. We have
$\MT(\tilde{s}_2)_{\ZZ_p}=k\MT(\tilde{s}_1)_{\ZZ_p}k^{-1}$. So we have
that all $\MT(\tilde{s})_{\ZZ_p}$ for $s$ in $\Sigma_p$ lie in one
$\GL(V_{\ZZ_p})$-orbit. 

The set of $\MT(\tilde{s})_{\FF_p}$ for $s$ in $\Sigma_p$ is contained
in one $\GL(V_{\FF_p})$-orbit, and hence is a finite set. If $s_1$ and
$s_2$ in $\Sigma_p$ are such that
$\MT(\tilde{s}_1)_{\FF_p}=\MT(\tilde{s}_2)_{\FF_p}$, then
$\MT(\tilde{s}_1)_{\ZZ_p}$ and $\MT(\tilde{s}_2)_{\ZZ_p}$ are
conjugated by an element of $G(\ZZ_p)$ by \cite[Exp.~XI,
Cor.~5.2]{Grothendieck1}, which says that the ``transporteur'' in
$G(\ZZ_p)$ between $\MT(\tilde{s}_1)_{\ZZ_p}$ and
$\MT(\tilde{s}_2)_{\ZZ_p}$ is smooth over~$\ZZ_p$.
\enddemo

Let $p$ be a prime with the following properties: 
\begin{itemize}
\item[1.] $p$ satisfies the conditions of Theorem~\ref{thm6.1};
\item[2.] $\Sigma_p$ is Zariski dense in~$Z$;
\item[3.] $G_{\FF_p}$ is smooth over~$\FF_p$;
\item[4.] $K_p=G_{\ZZ}(\ZZ_p)$;
\item[5.] All the $\MT(\tilde{s})_{\QQ_p}$ for $s$ in $\Sigma$ are split.
\end{itemize}
Indeed, the first four conditions only exclude finitely many primes,
and the last condition is verified by a set of primes of positive
density (Chebotarev). We replace $\Sigma$ by a suitable Zariski dense
subset of~$\Sigma_p$ such that the $\MT(\tilde{s})_{\ZZ_p}$ lie in one
$G_{\ZZ}(\ZZ_p)$-conjugacy class.

Let $s_0$ be an element of~$\Sigma$, and let $M=\MT(\tilde{s}_0)$.
The reciprocity morphism $r$ from $\Res_{F/\QQ}\Gm{F}{}$ to $M_\QQ$ is
surjective; hence $r((\QQ_p\otimes F)^*)$ is of finite index, say~$e$,
in~$M(\QQ_p)$. For each~$i$, the image $M_i$ of $M$ in $G_i$ is
nontrivial by the axioms for what constitutes a Shimura datum
(condition~(2.1.1.3) in~\cite[2.1]{De2}). By hypothesis, $M_{\QQ_p}$
and hence the $M_{i,\QQ_p}$ are split tori. For a split torus $T$
over~$\QQ_p$, $T(\QQ_p)$ modulo its maximal compact subgroup is a free
$\ZZ$-module of rank the dimension of~$T$. Since $M(\QQ_p)$ modulo its
maximal compact subgroup is not equal to a finite union of proper sub
$\ZZ$-modules, we can take an element $m$ of $M(\QQ_p)$ that satisfies
the following conditions:
\begin{itemize}
\item[1.] $m$ is in the image of multiplication by $e$ on~$M(\QQ_p)$;
\item[2.] For every simple factor $G_i$ of~$G$, the image of $m$ in
$G_i(\QQ_p)$ is not contained in a compact subgroup.
\end{itemize}
We will show that $Z\subset T_mZ$, so that $Z$ is of Hodge type by
Theorem~\ref{thm6.1}. Let $s$ be in~$\Sigma$. We have:
$$
T_ms=\{\ol{(\tilde{s},k_1mk_2)}\;|\;k_1\in K_p,\;k_2\in K_p\}. 
$$
By hypothesis, there exists an $h$ in $G_{\ZZ}(\ZZ_p)$ such that
$\MT(\tilde{s})_{\ZZ_p}=hM_{\ZZ_p}h^{-1}$. It follows that $hmh^{-1}$
is in $\MT(\tilde{s})(\QQ_p)$, and even in the image of the
reciprocity morphism for~$\tilde{s}$, so that: 
$$
\hbox{$\ol{(\tilde{s},hmh^{-1})}$ is in $\Gal(\Qbar/F){\cdot}s$.}
$$
We conclude that $\ol{(\tilde{s},hmh^{-1})}$ is in the intersection of
$T_ms$ and~$\Gal(\Qbar/F){\cdot}s$. As $T_ms$ is contained in~$T_mZ$,
and $\Gal(\Qbar/F){\cdot}s$ is contained in~$Z$,
$\ol{(\tilde{s},hmh^{-1})}$ is in $Z\cap T_mZ$. As both $Z$ and $T_mZ$
are defined over~$F$, we have: 
$$
\Gal(\Qbar/F){\cdot}s \subset Z\cap T_mZ.
$$
In particular, $s$ is in~$T_mZ$. Since this is so for all $s$
in~$\Sigma$, and $\Sigma$ is dense in~$Z$, we deduce that $Z$ is
contained in~$T_mZ$. This ends the proof of Theorem~1.2 in
the case where $i$ is bounded.

\vglue9pt {\it Remark} 7.3.2.
The proof of Theorem~1.2  just given, in the case
where the function $i$ is bounded, does not use the fact that $Z$ is a
curve. Hence it proves Theorem~1.2 without the condition that
$Z$ is a curve, but with the extra condition that $i$ is bounded. 

This kind of result can be useful. For example Vatsal and Cornut
(see~\cite{Vatsal1} and~\cite{Cornut1}) consider sets of Heegner
points on products of modular curves where the discriminants of the
endomorphism rings of the elliptic curves in question are products of
a fixed finite set of prime numbers (of course, the Andr\'e-Oort
conjecture in this case was already proved by Moonen
in~\cite[\S5]{Moonen3}, since it concerns a Shimura datum related to
moduli of abelian varieties).

\demo{{\rm 7.4.} The case where $i$ is not bounded}
In this section we assume that $i(s)$ is not bounded when $s$ ranges
through the points of~$\Sigma$. The strategy for proving
Theorem~1.2 in this case is as follows. For $s$ in
$\Sigma$ with $i(s)$ big enough, we show that there exist a prime
number $p$ and an element $m$ in $G(\QQ_p)$ that satisfy the first two
conditions in Theorem~\ref{thm6.1}, such that $Z\cap T_mZ$ contains
$\Gal(\Qbar/F){\cdot}s$ and such that $|\Gal(\Qbar/F){\cdot}s|$
exceeds the ``intersection number'' of $Z$ and~$T_mZ$. Then it follows
that $Z$ and $T_mZ$ do not intersect properly, and hence that $Z$ is
contained in~$T_mZ$ (it is here that we use the fact that $Z$ is a
curve). Theorem~\ref{thm6.1} then says that $Z$ is of Hodge type. 

We cite the following result (see \cite[Thm.~7.2]{Edixhoven1}) that
bounds the intersection of $Z$ and its images under Hecke
correspondences, if finite.

\specialnumber{7.4.1}\proclaim{Theorem}\label{thm6.4.1}
Let $(G,X)$ be a Shimura datum{\rm ,} let $K_1$ and $K_2$ be compact open
subgroups of $G(\AAf)${\rm ,} and let $Z_1$ and $Z_2$ be closed subvarieties
of the Shimura varieties $S_1:=\Sh_{K_1}(G,X)_\CC$ and
$S_2:=\Sh_{K_2}(G,X)_\CC${\rm ,} respectively. Suppose that $Z_1$ or $Z_2$
is of dimension at most one. Then there exists an integer
$c$ such that for all $g$ in $G(\AAf)$ for which $T_gZ_1\cap Z_2$ is
finite{\rm ,} 
$$
|T_gZ_1\cap Z_2| \leq c\,\deg(\pi_1\colon S_g\to S_1), 
$$
where $S_g=\Sh_{K_g}(G,X)_\CC$ with $K_g=K_1\cap gK_2g^{-1}${\rm ,} and with
$T_g$ and $\pi_1$ the morphism that is induced by the inclusion of
$K_g$ in~$K_1$.
\endproclaim

Applied in our situation, this theorem gives the following result:

\specialnumber{7.4.2} \proclaim{{C}orollary} \label{cor6.4.2}
There exists an integer $c$ such that for all $m$ in $G(\AAf)$ with
$Z\cap T_mZ$ finite{\rm ,}
$$
|Z\cap T_mZ| < c\; |K/K\cap mKm^{-1}|.
$$
\endproclaim

We note that if $p$ is a prime such that $K_p=G(\ZZ_p)$, and if $m$ is
in~$G(\QQ_p)$, then $K/K\cap mKm^{-1}$ is the $K_p$-orbit of the
lattice~$m V_{\ZZ_p}$ in the set of $\ZZ_p$-lattices in~$V_{\QQ_p}$,
and hence $|K/K\cap mKm^{-1}|$ is at most
$|\GL(V_{\ZZ_p}){\cdot}mV_{\ZZ_p}|$. In order to get prime numbers $p$
and elements $m$ in $G(\QQ_p)$ that verify the first two conditions in
Theorem~\ref{thm6.1} and are such that $|K/K\cap mKm^{-1}|$ is not too
big, we prove the following.

\vglue4pt {\elevensc Proposition 7.4.3}.
{\it Let $n${\rm ,} $r${\rm ,} $e$ and $B$ be nonnegative integers. There exists an
integer $k$ with the following property. Let $p$ be a prime number and
let $M$ be a split subtorus of~$\GL_{n,\QQ_p}${\rm ,} such that $M_{\ZZ_p}$
{\rm (}\/obtained from Zariski closure in~$\GL_{n,\ZZ_p}${\rm )} is a torus. Suppose
that with respect to a suitable $\ZZ$\/{\rm -}\/basis of the character group
$X^*(M_{\QQ_p})$ all coordinates of the differences of the characters
that intervene in $\QQ_p^n$ have absolute value at most~$B$. Let
$q_i\colon M\to M_i$ be $r$ quotients of~$M${\rm ,} with each $M_i$
nontrivial. Then there exists an element $m$ in $M(\QQ_p)$ such that
no $q_i(m)$ lies in a compact subgroup of~$M_i(\QQ_p)${\rm ,} and such that
$|\GL(V_{\ZZ_p}){\cdot}m^eV_{\ZZ_p}|<p^k$.}
\vglue4pt

{\it Proof}.
Let $n$, $r$, $e$ and $B$ be given. Let $p$ be a prime number, let $M$
be a split subtorus of~$\GL_{n,\QQ_p}$, such that $M_{\ZZ_p}$ is a
torus. Then we have a direct sum decomposition into character spaces: 
\vglue4pt
\centerline{${\displaystyle 
\ZZ_p^n = \oplus_{\chi\in X^*(M)} L_\chi.
}$}
\vglue4pt
\noindent 
Let $d$ be the dimension of~$M$, and let $f\colon \ZZ^d\to X^*(M)$ be
an isomorphism such that for each pair $(\chi_1,\chi_2)$ with
$L_{\chi_1}\neq 0 \neq L_{\chi_2}$ one has
$\|f^{-1}(\chi_1-\chi_2)\|\leq B$, where $\|{\cdot}\|$
denotes the maximum norm on~$\ZZ^d$. We identify $X^*(M)$ and its dual
$X_*(M)$ with $\ZZ^d$ via $f$ and its dual. Then the kernels of the
$X_*(q_i)$ give us $r$ subgroups $S_i$ of~$\ZZ^d$ of rank less
than~$d$. For $T$ a split torus over $\QQ_p$ one has
$T(\QQ_p)=X_*(T)\otimes \QQ_p^*$, and hence the valuation map
$v_p\colon\QQ_p^*\to \ZZ$ induces an isomorphism from $T(\QQ_p)$
modulo its maximal compact subgroup to~$X_*(T)$. It follows that for
an element $m$ of $M(\QQ_p)$ no $q_i(m)$ is in a compact subgroup if
and only if the image of $m$ in $\ZZ^d$ avoids all~$S_i$.

The quotient $\ZZ^d/(r+2)\ZZ^d$ has $(r+2)^d$ elements, whereas the
union of the images of the $S_i$ has at most $(r+2)^{d-1}r$ elements,
which is less than~$(r+2)^d$. Let $x$ be an element of $\ZZ^d$ with
$\|x\|< r+2$ whose image in $\ZZ^d/(r+2)\ZZ^d$ is not in the
union of the images of the~$S_i$. Let $m$ be an element in $M(\QQ_p)$
with image $x$ in~$\ZZ^d$. Let $a$ and $b$ be the smallest and largest
integers such that $p^a\ZZ_p^n\subset m^e\ZZ_p^n\subset
p^b\ZZ_p^n$. Then the $\GL_n(\ZZ_p)$-orbit of $m^e\ZZ_p^n$ is contained
in the set of $\ZZ_p$-lattices between $p^a\ZZ_p^n$
and~$p^b\ZZ_p^n$. The number of such lattices is the number of
subgroups of $(\ZZ/p^{a-b}\ZZ)^d$, hence is at most $p^{(a-b)d^2}$
(use that every subgroup is generated    by $d$ elements).
\pagegoal=50pc

It remains to bound~$a-b$. On $L_\chi$, $m^e$ acts, up to an element
of~$\ZZ_p^*$, by multiplication by $p^{\ld f^{-1}\chi,ex\rd}$, where
$\ld{\cdot},{\cdot}\rd$ denotes the standard pairing. Hence $a$ and
$b$ are the maximum and minimum of the $\ld f^{-1}\chi,ex\rd$, with
$L_\chi\neq0$. It follows that $a-b\leq edB(r+1)$. As $d\leq n$, we
can take $k=en^3B(r+1)$.
\hfill\qed\eject

Applied to our situation, we get the following consequence (note that
the $\MT(\tilde{s})_\QQ$ with $s$ in $\Sigma$ lie in one
$\GL(V_\QQ)$-conjugacy class).

\specialnumber{7.4.4} \proclaim{{C}orollary} \label{cor6.4.4}
There exists an integer $k$ with the following property. Let $s$ be
in~$\Sigma${\rm ,} and let $p$ be a prime such that $\MT(\tilde{s})_{\QQ_p}$
is split{\rm ,} and $\MT(\tilde{s})_{\FF_p}$ is a torus. Then there exists
an element $m$ in~$\MT(\tilde{s})(\QQ_p)$ such that\/{\rm : }\/
\begin{itemize}
\ritem{1.} $m$ is in the image of the reciprocity morphism
$r_{\tilde{s}}\colon(\QQ_p\otimes F)^*\to\MT(\tilde{s})(\QQ_p)${\rm ;}
\ritem{2.}  For all~$i$, the image of $m$ in $G_i(\QQ_p)$ is not in a
compact subgroup\/{\rm ;}\/
\ritem{3.} If $Z\cap T_mZ$ is finite{\rm ,} then $|Z\cap T_mZ|\leq p^k$. 
\end{itemize}
\endproclaim
\pagegoal=48pc

Theorem~\ref{thm3.1} provides us with positive reals numbers $c_1$ and
$c_2$ such that for all $s$ in $\Sigma$,
$$
|\Gal(\Qbar/F){\cdot}s| > c_1c_2^{i(s)}i(s)!. 
$$
For $x$ in~$\RR$, let $\pi_{F,1}(x)$ be the number of prime numbers
$p\leq x$ such that $F$ is split over~$\QQ_p$.  Chebotarev's density
theorem (see \cite[Ch.~VIII, \S4]{Lang1}) says that for $x$ in $\RR$
large enough, one has:
$$
\pi_{F,1}(x)\geq \frac{1}{2[F:\QQ]}{\cdot}\frac{x}{\log(x)}.
$$
Elementary calculus shows that for all $y$ in $\ZZ$ large enough,
there exists $x$ in $\RR$ such that:
$$
\frac{1}{2[F:\QQ]}{\cdot}\frac{x}{\log(x)} > y\quad\hbox{and}\quad
x^k<c_1c_2^yy!.
$$
As the function $i\colon\Sigma\to\ZZ$ is not bounded, we conclude that
there exists an $s$ in~$\Sigma$, a prime number~$p$ and an element $m$
in $\MT(\tilde{s})(\QQ_p)$ such that $m$ and $p$ satisfy the first two
conditions of Theorem~\ref{thm6.1}, and even the last one because
$Z\cap T_mZ$ contains~$\Gal(\Qbar/F){\cdot}s$.

\demo{Acknowledgements}
We would like to thank Paula Cohen for pointing out to us that it
would be interesting to prove Conjecture~\ref{AndreOortconjecture} in
the case of a set of special points that is contained in one Hecke
orbit. We thank Paula Cohen and Gisbert W\"ustholz for making
preliminary versions of their preprint \cite{CohenWustholz1} available
to us, and for useful conversations and correspondence. We thank
J\"urgen Wolfart for explaining to us the problem in his
article~\cite{Wolfart1}. We thank Yves Andr\'e and Richard Pink for
the work they have done as referees of the thesis of the second
author. We thank Rutger Noot for useful discussions on Shimura
varieties. 
\enddemo

\end{document}